\documentclass[12pt]{article}
\usepackage{amssymb,amsfonts, amsmath,amsthm,epsfig,tabularx,yfonts}
\usepackage{caption}
\usepackage{float}
\restylefloat{table}
\usepackage{xcolor}

\usepackage{hyperref}  

\usepackage{caption} 
\usepackage{subcaption} 

\usepackage{url,epsf,graphicx,tikz}
\usetikzlibrary{decorations.markings}

\newtheorem{theorem}{Theorem}
\newtheorem{proposition}[theorem]{Proposition}

\newtheorem{observation}[theorem]{Observation}

\newtheorem{conjecture}[theorem]{Conjecture}
\newtheorem{question}[theorem]{Question}
\newtheorem{problem}[theorem]{Problem}

\newcommand{\dC}{{\overrightarrow{C}}}         

\newcommand{\Sz}{{\rm Sz}}

\title{Selected topics on Wiener index} 

\author
{ Martin Knor\thanks{Slovak University of Technology in Bratislava,
Faculty of Civil Engineering, Department of Mathematics, Bratislava,
Slovakia. E-Mail: \texttt{knor@math.sk}},\quad Riste
\v{S}krekovski\thanks{FMF, University of Ljubljana \& Faculty of Information Studies, Novo mesto \&
Institute of Mathematics, Physics and Mechanics, Ljubljana \&
University of Primorska, FAMNIT, Koper, Slovenia. E-Mail:
\texttt{skrekovski@gmail.com}},\quad Aleksandra Tepeh\thanks{Faculty
of Information Studies, Novo mesto \& Faculty of Electrical
Engineering and Computer Science, University of Maribor, Slovenia.
E-Mail: \texttt{aleksandra.tepeh@gmail.com}} }

\begin{document}

\tikzset{middlearrow/.style={
        decoration={markings,
            mark= at position 0.5 with {\arrow{#1}} ,
        },
        postaction={decorate}
    }
}

\tikzset{My Style/.style={draw, circle, fill=red, scale=0.4}}

\maketitle

{\abstract
{The Wiener index is defined as the sum of distances between all unordered pairs of vertices in a graph. It is one of the most recognized and well-researched topological indices, which is on the other hand still a very active area of research. This work presents a natural continuation of the paper \textsl{Mathematical aspects of Wiener index} (Ars Math.
Contemp., 2016) in which several interesting open questions on the topic were outlined. Here we collect answers gathered so far, give further insights on the topic of 
extremal values of Wiener index in different settings, and present further intriguing problems and conjectures.
}}

\bigskip
\bigskip


\medskip

 \textbf{Math. Subj. Class.:} 05C05, 05C12, 05C20, 05C92, 92E10

\section{Introduction}

The \emph{Wiener index}, $W(G)$, is a topological index of a connected
graph, defined as the sum of the lengths of the shortest paths between all
unordered pairs of vertices in the graph. In other words, for a connected graph
$$
W(G)=\displaystyle\sum_{\{u,v\}\in V(G)} d(u,v),
$$
where $d(u,v)$ denotes the distance between vertices $u$ and $v$ in $G$.
This graph invariant has been investigated by numerous
authors (see e.g.~\cite{surv1,petra,DM3,KS_chapter,mathasp,Gut-sur}) under a variety of other names like transmission, total status, sum of all distances, path number and Wiener  number of a graph. 
Due to its basic character and applicability, it has arisen in diverse contexts, including efficiency of information, sociometry, mass transport, cryptography, theory of communication, molecular structure, complex network topology and many more.

The index was originally introduced in 1947 by Harold Wiener for the purpose of
determining the approximation formula of the boiling point of paraffin \cite{Wiener}.
The definition of Wiener index in terms of distances between vertices of a
graph was first given by Hosoya \cite{Hosoya}.

The \textsl{transmission} (also called the \textsl{distance}) of $u\in G$ is 
$t_G(u) =\sum_{v\in  V(G)} d_G(u, v)$.
 Thus the Wiener index can be expressed as 
$$
W(G)=\frac{1}{2}\displaystyle\sum_{v\in V(G)} t_G(v).
$$
Another view on the Wiener index was presented in \cite{win-dimension} as follows.
Suppose that $\{t_G(u)\,|\,u\in V(G\}=\{d_1, d_2, \ldots, d_k\}$. Assume in addition that $G$ contains $t_i$ vertices whose transmission is $d_i$, $1\leq i\leq k$. Then the Wiener index of $G$ can be expressed as 
$$
W(G) =\frac{1}{2}\displaystyle\sum_{i=1}^k t_i d_i.
$$
We therefore say that the \textit{Wiener dimension} ${\rm dim}_W(G)$ of $G$ is $k$. That is, the Wiener dimension of a graph is the number of different transmissions of its vertices.

Fundamental properties regarding extremal values of Wiener index are already a part of the folklore. In \cite{n1} and later in many subsequent papers (e.g. \cite{n2,n3}) it was shown that for trees on $n$ vertices, the maximum
Wiener index is obtained for the path $P_n$, and the minimum
for the star $S_n$.
Thus, for every tree $T$ on $n$ vertices, it holds
$$
(n-1)^2 = W(S_n) \le W(T) \le W(P_n) = \binom{n+1}{3}\,.
$$
Since the distance between any two distinct vertices is at least one,
we have that among all graphs on $n$ vertices $K_n$ has the smallest
Wiener index. In general, removing (resp.~adding) of an edge from a connected graph results in increased (resp.~decreased) Wiener index, which leads to the observation that Wiener index of a connected graph is less than or equal to the Wiener index of its spanning tree.
Therefore, for any connected graph $G$ on $n$ vertices, it holds
$$
\binom{n}{2} = W(K_n) \le W(G) \le W(P_n)=\binom{n+1}{3}\,.
$$

Despite extensive literature on the Wiener index, many interesting and basic questions remain open. In our previous survey \cite{mathasp} we have exposed some of them that mainly pertain to extremal values of Wiener index in different settings. In this paper we continue with summarizing knowledge accumulated since then, and integrate some new conjectures, problems and
ideas for possible future work. 



\section{Minimum Wiener index for chemical graphs}

The {\em degree} $\deg_G(v)$ of a vertex $v \in V(G)$ in a graph $G$ is $|N_G(v)|$, where $N_G(v)$
denotes the neighborhood of $v$ in $G$.
The {\em maximum degree} of a graph $G$, 
$\max_{v\in V(G)}\deg_G(v)$, is denoted by $\Delta(G)$, and the {\em minimum degree}, $\min_{v\in V(G)}\deg_G(v)$, is denoted by $\delta(G)$.

\medskip

Since every atom has a certain valency, chemists are often interested in graphs with restricted degrees, which correspond to valencies. Particularly interesting is the class of 
\textit{chemical graphs}, i.e.~graphs for which the degrees of its vertices do not exceed $4$. In \cite{mi-chem} the authors addressed an ``overlooked'' problem of determining
the minimum value of Wiener index and corresponding extremal graphs among chemical graphs with prescribed number of vertices. Note that the upper bound for this class of graphs is attained by paths.

\begin{problem}\label{min-chemical}
Find all the chemical graphs $G$ on $n$ vertices with the minimum value of
Wiener index.
\end{problem}

\begin{figure}[h]
	\begin{center}
		\begin{tikzpicture}[scale=0.9,style=thick]
		\node [My Style, name=h]   at (-10.5,0.5) {};		
		\node [My Style, name=i]   at (-9.5,2) {};
		\node [My Style, name=j]   at (-7.5,2) {};		
		\node [My Style, name=k]   at (-6.5,0.5) {};
		\node [My Style, name=l]   at (-7.5,-1) {};
		\node [My Style, name=m]   at (-9.5,-1) {};
	
		\draw[](h)--(i)-- (j)--(k)--(l)--(m)--(h)--(j)--(l)--(h);
		\draw[](i)--(k)-- (m)--(i);

		\node [My Style, name=g]   at (-5,0) {};		
		\node [My Style, name=a]   at (-4.75,1.5) {};
		\node [My Style, name=c]   at (-1.75,1.5) {};		
		\node [My Style, name=d]   at (-1.5,0) {};
		\node [My Style, name=e]   at (-2.5,-1) {};
		\node [My Style, name=f]   at (-4,-1) {};
	  \node [My Style, name=b]   at (-3.25,2) {};
	
		
		\draw[](a)--(b)-- (c)--(d)--(e)--(f)--(g)--(a)--(c)--(e)--(g)--(b)--(d)--(f)--(a);	

		\node [My Style, name=o]   at (0.25,1.5) {};		
		\node [My Style, name=u]   at (1.75,2) {};
		\node [My Style, name=p]   at (3.25,1.5) {};		
		\node [My Style, name=r]   at (3.5,0) {};
		\node [My Style, name=s]   at (2.5,-1) {};
		\node [My Style, name=t]   at (1,-1) {};
		\node [My Style, name=n]   at (0,0) {};
	
		
		\draw[](n)--(o)--(u)--(p)--(r)--(s)--(t)--(n)--(p)--(o)--(r)--(t)--(u)--(s)--(n);
		\draw[](o)--(p);		
		\end{tikzpicture}
	\end{center}
	\caption{Extremal graphs for $n=6$ and $n=7$.}
	\label{n-6-7}
\end{figure}
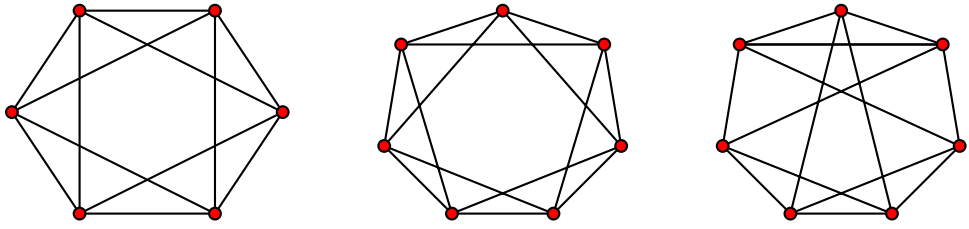

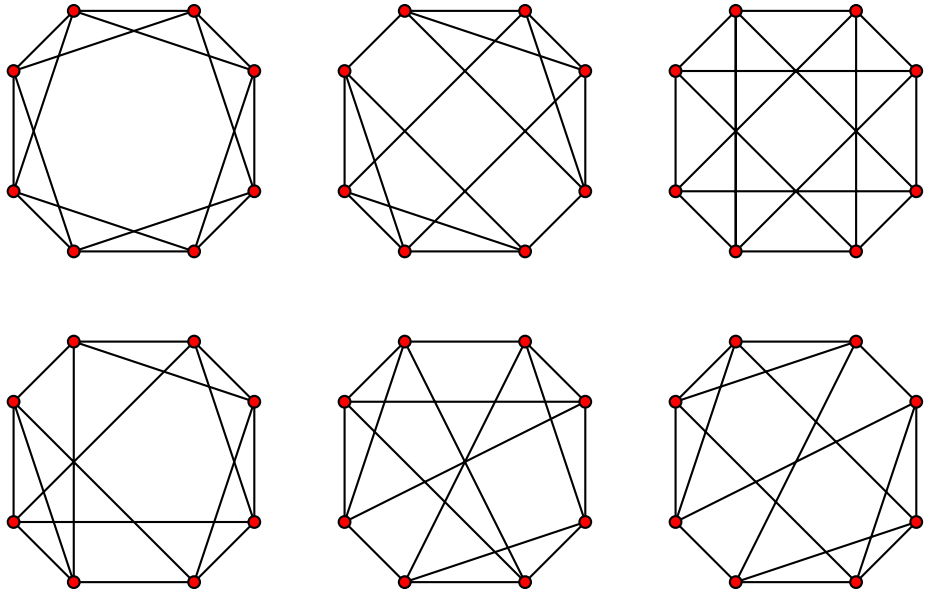
\begin{figure}[h]
	\begin{center}
		\begin{tikzpicture}[scale=0.8,style=thick]
		\node [My Style, name=a1]   at (1,4) {};		
		\node [My Style, name=a2]   at (3,0) {};
		\node [My Style, name=a3]   at (0,3) {};		
		\node [My Style, name=a4]   at (4,3) {};
		\node [My Style, name=a5]   at (1,0) {};
		\node [My Style, name=a6]   at (3,4) {};
		\node [My Style, name=a7]   at (0,1) {};
		\node [My Style, name=a0]   at (4,1) {};

		\draw[](a1)-- (a6)--(a4)--(a0)--(a2)--(a5)--(a7)--(a3)--(a1)--(a4)--(a2)--(a7)--(a1);
		\draw[](a3)--(a6)-- (a0)--(a5)--(a3);

		\node [My Style, name=b1]   at (9.5,3) {};		
		\node [My Style, name=b2]   at (5.5,3) {};
		\node [My Style, name=b3]   at (8.5,4) {};		
		\node [My Style, name=b4]   at (6.5,0) {};
		\node [My Style, name=b5]   at (9.5,1) {};
		\node [My Style, name=b6]   at (5.5,1) {};
		\node [My Style, name=b7]   at (6.5,4) {};
		\node [My Style, name=b0]   at (8.5,0) {};

		\draw[](b0)--(b4)--(b6)--(b2)--(b7)--(b3)--(b1)--(b5)--(b0)--(b2)--(b4)--(b1)--(b7)--(b5);
		\draw[](b3)--(b6);	
    \draw[](b3)--(b5);
		\draw[](b0)--(b6);
		
		\node [My Style, name=c1]   at (11,1) {};		
		\node [My Style, name=c2]   at (15,3) {};
		\node [My Style, name=c3]   at (12,4) {};		
		\node [My Style, name=c4]   at (12,0) {};
		\node [My Style, name=c5]   at (11,3) {};
		\node [My Style, name=c6]   at (15,1) {};
		\node [My Style, name=c7]   at (14,4) {};
		\node [My Style, name=c0]   at (14,0) {};

		\draw[](c4)--(c1)--(c5)--(c3)--(c7)--(c2)--(c6)--(c0)--(c4)--(c3)--(c6)--(c1)--(c7)--(c0)--(c5)--(c2)--(c4)--(c3);

		\node [My Style, name=d1]   at (0,-4.5) {};		
		\node [My Style, name=d2]   at (3,-5.5) {};
		\node [My Style, name=d3]   at (3,-1.5) {};		
		\node [My Style, name=d4]   at (4,-2.5) {};
		\node [My Style, name=d5]   at (1,-5.5) {};
		\node [My Style, name=d6]   at (4,-4.5) {};
		\node [My Style, name=d7]   at (0,-2.5) {};
		\node [My Style, name=d0]   at (1,-1.5) {};

		\draw[](d0)--(d3)--(d4)--(d6)--(d2)--(d5)--(d1)--(d7)--(d0)--(d4)--(d2)--(d7)--(d5)--(d0);
		\draw[](d3)--(d6)--(d1)--(d3);

		\node [My Style, name=e1]   at (9.5,-2.5) {};		
		\node [My Style, name=e2]   at (6.5,-5.5) {};
		\node [My Style, name=e3]   at (5.5,-2.5) {};		
		\node [My Style, name=e4]   at (8.5,-5.5) {};
		\node [My Style, name=e5]   at (8.5,-1.5) {};
		\node [My Style, name=e6]   at (5.5,-4.5) {};
		\node [My Style, name=e7]   at (9.5,-4.5) {};
		\node [My Style, name=e0]   at (6.5,-1.5) {};

    \draw[](e0)--(e5)--(e1)--(e7)--(e4)--(e2)--(e6)--(e3)--(e0)--(e4)--(e3)--(e1)--(e6)--(e0);
		\draw[](e2)--(e5)--(e7)--(e2);	
		
		\node [My Style, name=f1]   at (12,-5.5) {};		
		\node [My Style, name=f2]   at (15,-2.5) {};
		\node [My Style, name=f3]   at (12,-1.5) {};		
		\node [My Style, name=f4]   at (14,-5.5) {};
		\node [My Style, name=f5]   at (14,-1.5) {};
		\node [My Style, name=f6]   at (11,-4.5) {};
		\node [My Style, name=f7]   at (15,-4.5) {};
		\node [My Style, name=f0]   at (11,-2.5) {};

    \draw[](f0)--(f3)--(f5)--(f2)--(f7)--(f4)--(f1)--(f6)--(f0)--(f5)--(f1)--(f7)--(f3)--(f6)--(f2)--(f4)--(f0);

		\end{tikzpicture}
	\end{center}
	\caption{Extremal graphs for $n=8$.}
	\label{n-8}
\end{figure}

Inserting of an edge in a graph decreases the Wiener index, thus one would expect that its minimum in the class of chemical graphs is attained by $4$-regular graphs. Using  a computer it was verified that for $n\in \{1,2,\ldots,5\}$ minimum is attained for $K_n$. Extremal graphs in cases $n=6,7$ are presented in Figure \ref{n-6-7}. Observe that the first two graphs in this figure are circulant graphs $C_6(1,2)$ and $C_7(1,2)$, respectively, and they are vertex-transitive.
There are $1929$ simple connected graphs on $8$ vertices and the minimum Wiener index
value is $40$, which is attained by only $6$ graphs depicted in Figure \ref{n-8}. Note that the first three graphs, which are the circulant graph $C_8(1,2)$, the Cartesian product $K_4 \Box P_2$ and the complete bipartite graph $K_{4,4}=C_8(1,3)$, respectively, are vertex-transitive. The above cases support the following conjecture.

\begin{conjecture}
Every chemical graph $G$ on $n\ge 5$ vertices with the minimum value of
Wiener index is 4-regular.
\end{conjecture}

Although computer results indicate the above conjecture to be true, the problem seems to be far from tractable. In \cite{mi-chem} it is shown that a chemical graph with the minimum value of Wiener index has at most $3$ vertices of degree smaller than $4$. In fact, a more general statement holds.

\begin{observation}
If $G$ is a graph on $n$ vertices with maximum degree $\Delta$, $n \geq  \Delta + 1$, and
with the minimum possible value of Wiener index,  then $G$ contains at most $\Delta - 1$ vertices whose degree is strictly smaller than $\Delta$, and these vertices induce a clique.
\end{observation}



\section{Prescribed degrees}

As mentioned earlier, among $n$-vertex graphs with minimum degree at least $1$,
the maximum Wiener index is attained by $P_n$.
But when restricting to minimum degree at least $2$, the extremal graph is different.
Observe that with the reasonable assumptions $\Delta \ge 2$ and $\delta \le n-1$, 
the following holds:
\begin{itemize}
\item $W(P_n)=\max\{W(G);\,
G\mbox{ has maximum degree at most }\Delta \mbox{ and }n\mbox{ vertices}\}$,
\item $W(K_n)=\min\{W(G);\, G\mbox{ has minimum degree at least }\delta \mbox{ and }n\mbox{ vertices}\}$.
\end{itemize}
Analogous reasons motivate the following two problems from \cite{mathasp}. 

\begin{problem}
\label{prob:up-bound}
What is the maximum Wiener index among $n$-vertex graphs with minimum degree
at least $\delta$?
\end{problem}

\begin{problem}
\label{prob:lo-bound}
What is the minimum Wiener index among $n$-vertex graphs with maximum degree
at most $\Delta$?
\end{problem}

Both problems are still on the list of unsolved problems, but several results were obtained under additional requirements.  
Fischermann et al.~\cite{fish}, and independently Jelen and Trisch~\cite{jel,tri} solved Problem \ref{prob:lo-bound} for trees.
In addition, they determined the trees which maximize the Wiener index among all trees of given order whose vertices are either endvertices or of maximum degree $\Delta$. 

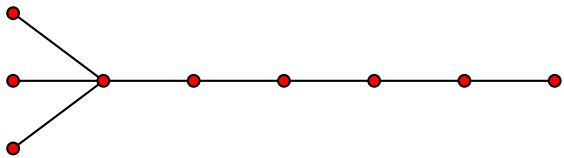
\begin{figure}[ht!]
\begin{center}
\begin{tikzpicture}[scale=1.2,style=thick]

		\node [My Style, name=a]   at (0,0) {};		
		\node [My Style, name=b]   at (1,0) {};
		\node [My Style, name=c]   at (2,0) {};		
		\node [My Style, name=d]   at (3,0) {};
		\node [My Style, name=e]   at (4,0) {};
		\node [My Style, name=f]   at (5,0) {};
		\node [My Style, name=g]   at (6,0) {};
		\node [My Style, name=i]   at (0,0.75) {};
		\node [My Style, name=j]   at (0,-0.75) {};
		
		\draw[](a)--(b)--(c)--(d)--(e)--(f)--(g);
		\draw[](i)--(b)--(j);

\end{tikzpicture}
\end{center}
\caption{Graph $T_{9,4}$.}

\label{stev}
\end{figure}

Stevanovi\'{c}~\cite{stev} solved Problem  \ref{prob:lo-bound} for trees under the assumption that the maximum degree is precisely $\Delta$. Let $T_{n,\Delta}$ be the tree on $n$ vertices obtained by taking a path on $n-\Delta+1$ vertices and joining new $\Delta-1$ vertices to one end-vertex of the path, see Figure \ref{stev}.

\begin{theorem}
For every $n$-vertex graph $G$ with maximum degree $\Delta \geq 2$ it holds that $W(G)\leq W(T_{n,\Delta})$ with equality if and only if $G$ is $T_{n,\Delta}$.
\end{theorem}


Dong and Zhou  \cite{DongZhou} determined the maximum Wiener index of unicyclic
graphs with fixed maximum degree and they characterized the unique extremal graph.

Lin  \cite{Lin14} characterized trees with the maximal Wiener index in the class of trees of order $n$ with exactly $k$ vertices of maximum degree, and proposed   analogous problem for the minimum. The solution of this problem was recently presented by Bo\v{z}ovi\'{c} et al.~in \cite{bozo}. The same authors considered a similar problem with a predetermined value of the maximum degree, i.e.~they obtained the maximal value of Wiener
index in the class of trees of order $n$ with exactly $k$ vertices of a given maximum degree and showed that the corresponding maximal
trees are caterpillars with certain properties.



Recently Alochukwu and Dankelmann \cite{AloDankel} obtained the following asymptotically sharp upper bound in terms of given minimum and maximum degree.

\begin{theorem}
Let $G$ be a graph of order $n$, minimum degree $\delta$ and maximum degree $\Delta$. Then $W(G) \leq \binom{n-\Delta+\delta}{2}\frac{n+2\Delta}{\delta+1}+2n(n-1)$,  and this bound is sharp apart from an additive constant.
\end{theorem}

Another interesting class of graphs with restrictions on degrees is 
the class of \textit{regular graphs}, i.e.~graphs for which $\Delta(G)=\delta(G)$.
In general, introducing edges in a graph decreases the
Wiener index, but in the class of $r$-regular graphs on $n$ vertices the number of edges is fixed, therefore the following conjecture from 
\cite{miAMC} seems to be reasonable. The {\em diameter}, ${\rm diam}(G)$, of a graph $G$ is the maximum distance between all pairs of vertices, i.e.~${\rm diam}(G)=\max \{d(u,v)|\,\, u,v \in V(G)\}$.

\begin{conjecture}
Among all $r$-regular graphs on $n$ vertices, the maximum Wiener index
is attained by a graph with the maximum possible diameter.
\end{conjecture}

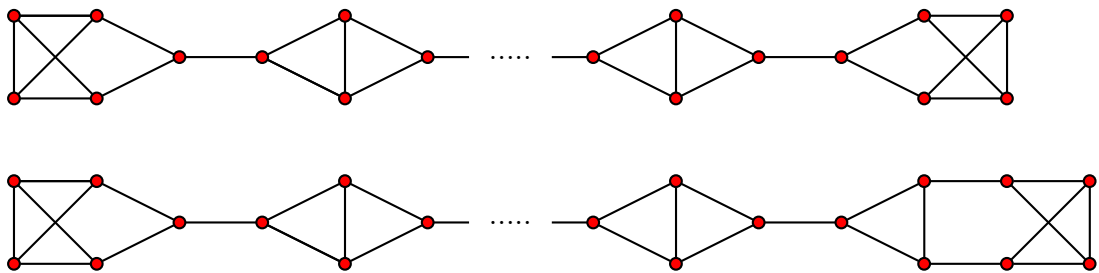
\begin{figure}[ht!]
\begin{center}
\begin{tikzpicture}[scale=1.1,style=thick]

		\node [My Style, name=a]   at (0,0) {};		
		\node [My Style, name=b]   at (1,0) {};
		\node [My Style, name=c]   at (2,0.5) {};		
		\node [My Style, name=d]   at (1,1) {};
		\node [My Style, name=e]   at (0,1) {};
		\node [My Style, name=f]   at (3,0.5) {};
		\node [My Style, name=g]   at (4,0) {};
		\node [My Style, name=h]   at (4,1) {};
		\node [My Style, name=i]   at (5,0.5) {};
		\node [My Style, name=j]   at (7,0.5) {};
		\node [My Style, name=k]   at (8,0) {};
		\node [My Style, name=l]   at (8,1) {};
		\node [My Style, name=m]   at (9,0.5) {};
		\node [My Style, name=n]   at (10,0.5) {};
		\node [My Style, name=o]   at (11,0) {};
		\node [My Style, name=p]   at (12,0) {};
		\node [My Style, name=r]   at (13,0) {};
		\node [My Style, name=s]   at (13,1) {};
		\node [My Style, name=t]   at (12,1) {};
		\node [My Style, name=u]   at (11,1) {};

		\draw[](c)--(d)--(e)--(a)--(b)--(c)--(f)--(g)--(i);
		\draw[](h)--(f)--(g)--(h)--(i);
		\draw[](a)--(d)--(e);
		\draw[](e)--(b);
		\draw[](j)--(l)--(k)--(j);
		\draw[](l)--(m)--(k);
		\draw[](m)--(n)--(o)--(p)--(r)--(s)--(t)--(u)--(n);
		\draw[](t)--(r);
		\draw[](p)--(s);
		\draw[](o)--(u);
		
\draw[](i)--(5.5,0.5);	
\draw (6,0.5) node {$.....$};	
\draw[](6.5,0.5)--(j);

		\node [My Style, name=a1]   at (0,2) {};		
		\node [My Style, name=b1]   at (1,2) {};
		\node [My Style, name=c1]   at (2,2.5) {};		
		\node [My Style, name=d1]   at (1,3) {};
		\node [My Style, name=e1]   at (0,3) {};
		\node [My Style, name=f1]   at (3,2.5) {};
		\node [My Style, name=g1]   at (4,2) {};
		\node [My Style, name=h1]   at (4,3) {};
		\node [My Style, name=i1]   at (5,2.5) {};
		\node [My Style, name=j1]   at (7,2.5) {};
		\node [My Style, name=k1]   at (8,2) {};
		\node [My Style, name=l1]   at (8,3) {};
		\node [My Style, name=m1]   at (9,2.5) {};
		\node [My Style, name=n1]   at (10,2.5) {};
		\node [My Style, name=o1]   at (11,2) {};
		\node [My Style, name=p1]   at (12,2) {};
		\node [My Style, name=t1]   at (12,3) {};
		\node [My Style, name=u1]   at (11,3) {};

		\draw[](c1)--(d1)--(e1)--(a1)--(b1)--(c1)--(f1)--(g1)--(i1);
		\draw[](h1)--(f1)--(g1)--(h1)--(i1);
		\draw[](a1)--(d1)--(e1);
		\draw[](e1)--(b1);
		\draw[](j1)--(l1)--(k1)--(j1);
		\draw[](l1)--(m1)--(k1);
		\draw[](m1)--(n1)--(o1)--(p1)--(t1)--(u1)--(n1);
		\draw[](t1)--(o1);
		\draw[](p1)--(u1);

		
\draw[](i1)--(5.5,2.5);	
\draw (6,2.5) node {$.....$};	
\draw[](6.5,2.5)--(j1);		
		
\end{tikzpicture}
\end{center}
\caption{Graphs $L_{4k+2}$ (above) and $L_{4k+4}$ (below).}

\label{ln}
\end{figure}

The above conjecture can be supported by the fact that in the case of
trees, where the number of edges is fixed as well, the maximum Wiener index is attained by $P_n$ which has the largest diameter. In fact, Chen et al.~\cite{chen1} recently proved that the conjecture is valid for $r=3$. More precisely, they proved a conjecture from \cite{miAMC}, that cubic graphs of the form $L_n$, presented in Figure~\ref{ln}, have maximum Wiener index among all cubic graphs of order $n$.


The minimum Wiener index in the class of trees is attained by $S_n$, which has the smallest diameter. A similar claim may hold for regular graphs \cite{miAMC}.

\begin{conjecture}
Among all $r$-regular graphs on $n$ vertices, the minimum Wiener index
is attained by a graph with the minimum possible diameter.
\end{conjecture}

Finally, the following problem from \cite{mi-chem} is of a special interest.

\begin{problem}
\label{reg}
Find all $k$-regular graphs on $n$ vertices with the smallest value of
Wiener index.
\end{problem}

As observed in \cite{mi-chem}, Problem \ref{reg}
is surprisingly related to the cages and the following famous degree-diameter problem
(see~\cite{MS} for details).

\begin{problem}[The degree-diameter problem]
\label{prob:chemic}
Given positive integers $d$ and $k$, find the largest possible number $n(d,k)$ of vertices in a graph of maximum degree $d$ and diameter $k$.
\end{problem}

Computer results in \cite{mi-chem} (see also \cite{dW}) showed that among graphs with the minimum Wiener index there are graphs achieving $n(k,d)$ for pairs $(k,d)$ from $\{(3,2),(3,3),(4,2)\}$. 
There might appear graphs achieving $n(k,d)$ also for higher values of
diameter $d$, but for those we could not search the space of $k$-regular
graphs of order $n$ exhaustively.
Anyway, for higher diameters the graphs achieving $n(k,d)$ do not need to be
those with the smallest Wiener index. Among extremal graphs found by a computer, $n(3,2)$ and $n(3,3)$ are realized by the well-known Petersen
graph and the Flower snark $J_5$.
Interestingly, there appears also the Heawood graph, which is the
Cage$(3,6)$, i.e., the smallest graph of degree $3$ and girth $6$,
see~\cite{EJ}.

\medskip

The following conjectures were proposed in \cite{mi-chem}
(probably, it suffices to choose $n_k=k+1$ therein).

\begin{conjecture}[The even case conjecture]
\label{conj:eregular}
Let $k\ge 3$, and let $n$ be large enough with respect to $k$,
say $n\ge n_k$.
Suppose that $G$ is a graph on $n$ vertices with the maximum degree $k$,
and with the smallest possible value of Wiener index.
If $kn$ is even, then $G$ is $k$-regular.  
\end{conjecture}

\begin{conjecture}[The odd case conjecture]
\label{conj:oregular}
Let $k\ge 3$, and let $n$ be large enough with respect to $k$,
say $n\ge n_k$.
Suppose that $G$ is a graph on $n$ vertices with the maximum degree $k$,
and with the smallest possible value of Wiener index.
If $kn$ is odd, then $G$ has a unique vertex of degree smaller than $k$
and in that case this smaller degree is $k-1$.
\end{conjecture}

%
%
\section{Wiener index of digraphs} 
\label{4}

A {\em directed graph} (a {\em digraph}) $D$ is given by a set of
vertices $V(D)$ and a set of ordered pairs of vertices $A(D)$ called
{\em directed edges} or {\em arcs}. If $uv$ is an arc in $D$, we say that $u$ \textit{dominates} $v$.
The \textit{out-degree} $d^+(u)$ of a vertex $u\in V(D)$ is the number of its \textit{out-neighbors}, i.e.~the vertices, dominated by $u$.
A {\em (directed) path} in $D$ is a sequence of vertices
$v_0,v_1,\ldots,v_k$ such that $v_{i-1}v_i$ is an arc of $D$ for each
$i\in \{1,2,\ldots, k-1\}$. The {\em distance} $d(u,v)$ between vertices $u,v\in V(D)$ is
the length of a shortest path from $u$ to $v$. Notice that $d(u,v)$ is usually distinct from $d(v,u)$. 

Early studies of Wiener index of digraphs were limited to \textit{strongly connected} digraphs, i.e.~digraphs in which a directed path between every pair of vertices exists.  
However, in the studies of real directed networks it is possible that there is no directed path connecting some pairs of vertices, thus the convention $d(u,v)=0$ is used if there is no directed path from $u$ to $v$ \cite{bonc1,bonc2}. Under this assumption, 
in analogy to graphs, the Wiener index $W(D)$ of a digraph $D$ is
defined as the sum of all distances, where each ordered
pair of vertices is taken into account.
Hence,
$$
W(D)=\sum_{(u,v)\in V(D)\times V(D)}d(u,v).
$$

Let $W_{\max}(G)$ and $W_{\min}(G)$ be the maximum possible and the minimum possible, respectively, Wiener index among all digraphs
obtained by orienting the edges of a graph $G$. If an orientation of $G$  achieves the minimum Wiener index $W_{\min}(G)$, we call this orientation a minimum Wiener index orientation of $G$.

\begin{problem}
For a given graph $G$ find  $W_{\max}(G)$ and $W_{\min}(G)$.
\end{problem}

In \cite{KST2} there was posed a question if it is NP-hard to find an orientation of a given graph which maximizes the Wiener index. Dankelmann \cite{dankelmann22} answered it affirmatively.
Plesn{\'\i}k \cite{P} proved that finding a strongly connected orientation of a given graph $G$ that minimizes the Wiener index is NP-hard too, but the case for non-necessarily strongly connected digraphs is unsolved \cite{KST2} in general. However, it can be decided in polynomial time if a given graph with $m$ edges has an orientation which Wiener index is precisely $m$ (note that it cannot be less). 

\begin{problem}
What is the complexity of finding
 $W_{\min}(G)$ for an input graph $G$?
\end{problem}

The following conjecture from \cite{KST2} remains unsolved as well, but it is known to hold for bipartite graphs, unicyclic graphs, the Petersen graph and prisms.

\begin{conjecture}
\label{conj:acyclic}
For every graph $G$, the value $W_{\min}(G)$ is achieved by some acyclic
orientation of $G$.
\end{conjecture}

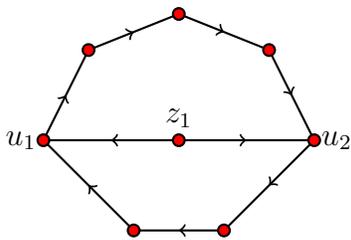
\begin{figure}[ht!]
\begin{center}
\begin{tikzpicture}[scale=0.6,style=thick]


\node [My Style, name=a1]   at (-2,0) {};
\node [My Style, name=a2]   at (0,0.8) {};
\node [My Style, name=a3] at (2,0) {};
\node [My Style, name=u1]   at (-3,-2) {};
\node [My Style, name=u2]   at (3,-2) {};
\node [My Style, name=b1] at (-1,-4) {};
\node [My Style, name=b2]   at (1,-4) {};
\node [My Style, name=z1]   at (0,-2) {};

\draw[middlearrow={>}] (u1) -- (a1); \draw[middlearrow={>}] (a1) -> (a2); \draw[middlearrow={>}] (a2) -> (a3);
\draw[middlearrow={>}] (a3) -> (u2); \draw[middlearrow={>}] (u2) -> (b2); \draw[middlearrow={>}] (b2) -> (b1);
\draw[middlearrow={>}] (b1) -> (u1); \draw[middlearrow={>}] (z1) -> (u1); \draw[middlearrow={>}] (z1) -> (u2);

\draw (3.5,-2) node {$u_2$};
\draw (-3.5,-2) node {$u_1$};
\draw (0,-1.5) node {$z_1$};

\end{tikzpicture}

\end{center}

\caption{An orientation of $\Theta_{3,2,1}$.}

\label{theta321}
\end{figure}

In \cite{moon,P} Plesn{\'\i}k and Moon found strongly connected tournaments (orientations of $K_n$) with the maximum and the second maximum Wiener index.
In \cite{KST} it was shown that the same tournaments solve the problem if we drop out the requirement that the digraph should be strongly connected. 
In the same paper oriented
$\Theta$-graphs are studied. By $\Theta_{a,b,c}$ we denote a graph obtained when
two distinct vertices are connected by three
internally vertex-disjoint paths of lengths $a+1$, $b+1$ and $c+1$,
respectively, where $a\geq b\geq c$ and $b\geq 1$ (see Figure \ref{theta321} where a non-strongly connected orientation of $\Theta_{3,2,1}$ is depicted).
Although intuitively one may expect that $W_{\max}$ is attained for some strongly connected orientation, this is not the case. Namely, in \cite{KST} it is shown that the orientation of $\Theta_{a,b,c}$ which achieves the maximum Wiener index is not 
strongly connected if $c\ge 1$. 

For strongly connected
orientations of $\Theta_{a,b,c}$, it was shown that the maximum Wiener index is achieved by the one in which the union of the $u_1,u_2$-paths
of lengths $a+1$ and $b+1$ forms a directed cycle.
Li and Wu \cite{LiWu} confirmed the conjecture from \cite{KST}, that the same holds if we drop the assumption that orientations are strongly connected.

\begin{theorem} 
Let $a\geq b \geq c$.
Then $W_{\max}(\Theta_{a,b,c})$ is attained by an orientation of
$\Theta_{a,b,c}$ in which the union of the paths of lengths $a+1$
and $b+1$ forms a directed cycle.
\end{theorem}

\hspace{0.5cm}

However, the following conjecture remains open. 

\begin{conjecture}
Let $G$ be a 2-connected chordal graph. Then $W_{\max}(G)$ is attained by an orientation
which is strongly connected.
\end{conjecture}

Among digraphs on $n$ vertices, the directed cycle $\dC_n$ (in which all edges are directed in the same way, say clockwise) achieves the maximum Wiener index. In \cite{KST3} digraphs with the second maximum Wiener index were investigated. In \cite{KST2} the Wiener theorem was generalized to directed graphs, as well as a relation between the Wiener index and betweenness centrality.

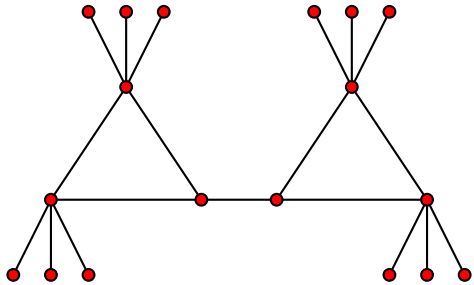
\begin{figure}[ht!]
\begin{center}
\begin{tikzpicture}[scale=1,style=thick]

		\node [My Style, name=a]   at (0,0) {};		
		\node [My Style, name=b]   at (2,0) {};
		\node [My Style, name=c]   at (3,0) {};		
		\node [My Style, name=d]   at (5,0) {};
		\node [My Style, name=e]   at (4,1.5) {};
		\node [My Style, name=f]   at (1,1.5) {};
		\node [My Style, name=g]   at (0.5,2.5) {};
		\node [My Style, name=h]   at (1,2.5) {};
		\node [My Style, name=i]   at (1.5,2.5) {};
		\node [My Style, name=j]   at (3.5,2.5) {};
		\node [My Style, name=k]   at (4,2.5) {};
		\node [My Style, name=l]   at (4.5,2.5) {};
		\node [My Style, name=m]   at (-0.5,-1) {};
		\node [My Style, name=n]   at (0,-1) {};
		\node [My Style, name=o]   at (0.5,-1) {};
		\node [My Style, name=p]   at (4.5,-1) {};
		\node [My Style, name=r]   at (5,-1) {};
		\node [My Style, name=s]   at (5.5,-1) {};

		\draw[](b)--(f)--(a)--(b)--(c)--(d)--(e)--(c);
		\draw[](f)--(i);
		\draw[](g)--(f)--(h);
		\draw[](e)--(l);
		\draw[](j)--(e)--(k);
		\draw[](a)--(o);
		\draw[](m)--(a)--(n);
		\draw[](d)--(s);
		\draw[](p)--(d)--(r);
		
\end{tikzpicture}
\end{center}
\caption{A graph $G$, for which $W_{\min}(G)$ is not achieved for any
$\chi(G)$-coloring-induced orientation of $G$.}

\label{counterex}
\end{figure}

An orientation of a graph $G$ is called \textit{$k$-coloring-induced}, if it is obtained from a proper $k$-coloring of $G$ such that each edge
is oriented from the end-vertex with the bigger color to the end-vertex with the smaller color. In \cite{KST2} it was proved that graphs with at most one cycle and prisms attain the minimum Wiener index for $k$-coloring-induced orientation with $k$ being the chromatic number $\chi(G)$. The same holds for bipartite graphs, complete graphs, Petersen graph and others. These observations lead to the conjecture that $W_{\min}(G)$ of an arbitrary graph is achieved for a $\chi(G)$-coloring-induced orientation, which Fang and Gao \cite{counterex} showed to be false. 
They expressed the Wiener index of a digraph $D$ as $W(D)=\sum_{u\in V(D)}w(u)$ where $w(u)=\sum_{v\in V(D)}d(u,v)$, and defined the notion of Wiener increment. For $u\in V(D)$ the \textit{Wiener increment of $u$} is defined as $\Delta w(u)=w(u)-d^+(u)$. The \textit{Wiener increment of $D$}, $\Delta W(D)$, is the sum of Wiener increments of all vertices of $D$. Fang and Gao observed that the comparison of Wiener indices of two different orientations of a graph is equal to the comparison of their Wiener increments. Using this observation they found that for the graph $G$ in Figure~\ref{counterex}, $W_{\min}(G)$ cannot be achieved for any $\chi(G)$-coloring-induced orientation of $G$, and this is not the only counterexample. Moreover, their investigations lead them to pose the following two conjectures.


\begin{conjecture}
For any given constant $k \geq 3$, there exists a $3$-colorable graph $G$ such that any minimum Wiener index orientation of $G$ has a directed path of length $k$.
\end{conjecture}

\begin{conjecture}
For any given constant $k \geq 3$, there exists a $3$-colorable graph $G$ such that $W_{\min}(G)$ cannot be achieved by
any $k$-coloring-induced orientation.
\end{conjecture}

\hspace{0.5cm}

\begin{figure}[ht!]
\begin{center}
\begin{tikzpicture}[scale=1.2,style=thick]


\node [My Style, name=a2]   at (-4,0) {};
\node [My Style, name=a1]   at (-3,0) {};
\node [My Style, name=a]   at (-2,0) {};
\node [My Style, name=b]   at (-1,0) {};
\node [My Style, name=c] at (0,0) {};
\node [My Style, name=d]   at (1,0) {};
\node [My Style, name=e]   at (2,0) {};

\node [My Style, name=f1]   at (3,0) {};
\node [My Style, name=f]   at (4,0) {};
\node [My Style, name=g]   at (5,0) {};
\node [My Style, name=h] at (6,0) {};
\node [My Style, name=i]   at (7,0) {};
\node [My Style, name=j]   at (8,0) {};
\node [My Style, name=j1]   at (9,0) {};

\draw[middlearrow={>}] (f) -> (g);
\draw[middlearrow={>}] (h) -> (g);
\draw[middlearrow={>}] (i) -> (h);
\draw[middlearrow={>}] (a1) -> (a); 
\draw[middlearrow={>}] (j) -> (i);
\draw[middlearrow={>}] (a2) -> (a1); 
\draw[middlearrow={>}] (a) -> (b); 
\draw[middlearrow={>}] (b) -> (c); 
\draw[middlearrow={>}] (c) -> (d);
\draw[middlearrow={>}] (d) -> (e); 
\draw[middlearrow={>}] (j) -> (j1);
\draw[middlearrow={>}] (f1) -> (f);

\end{tikzpicture}
\end{center}
\caption{A no-zig-zag path (left) and a zig-zag path (right) on six vertices.}

\label{paths}
\end{figure}
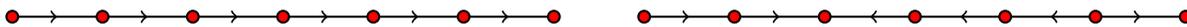

\begin{figure}[ht!]
\begin{center}
\begin{tikzpicture}[scale=0.8,style=thick]


\node [My Style, name=a]   at (-1,1) {};
\node [My Style, name=b]   at (1,1) {};
\node [My Style, name=c] at (0,0) {};
\node [My Style, name=d]   at (0,-1.5) {};
\node [My Style, name=e]   at (-1,-2.5) {};
\node [My Style, name=f] at (1,-2.5) {};

\draw[middlearrow={>}] (a) -- (c); \draw[middlearrow={>}] (b) -> (c); \draw[middlearrow={>}] (c) -> (d);
\draw[middlearrow={>}] (d) -> (e); \draw[middlearrow={>}] (d) -> (f); 


\node [My Style, name=a1]   at (4,1) {};
\node [My Style, name=b1]   at (6,1) {};
\node [My Style, name=c1] at (5,0) {};
\node [My Style, name=d1]   at (5,-1.5) {};
\node [My Style, name=e1]   at (4,-2.5) {};
\node [My Style, name=f1] at (6,-2.5) {};

\draw[middlearrow={>}] (a1) -- (c1); \draw[middlearrow={>}] (b1) -> (c1); \draw[middlearrow={<}] (c1) -> (d1);
\draw[middlearrow={>}] (d1) -> (e1); \draw[middlearrow={>}] (d1) -> (f1); 


\end{tikzpicture}
\end{center}
\caption{The graph on the left-hand side has
two core vertices, while the right-hand side one has no core vertex.}

\label{cores}
\end{figure}
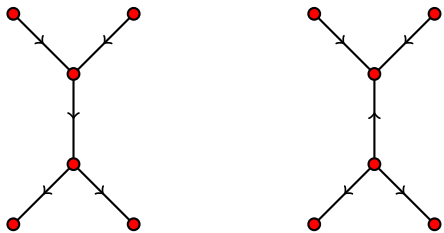

In \cite{KST2} orientations of trees with the maximum Wiener index were considered. An orientation of a tree is called {\em no-zig-zag} if there is no subpath in which edges change the orientation twice, see Figure \ref{paths}. 
 A different view on no-zig-zag trees can be described as follows. A vertex $v$ in a directed tree $T$ is {\em core}, if for every vertex $u$ of $T$ there exists either a directed path from $u$ to
$v$ or a directed path from $v$ to $u$, see Figure \ref{cores}. 
Notice that then in each component $C$ of
$T-v$ all edges point in the direction towards $v$ or all edges
point in the direction from $v$.

\hspace{0.3cm}

In \cite{KST2} the following conjecture was proposed

\begin{conjecture} 
\label{conj:zig-zag} Let $T$ be a tree. Then every orientation of $T$ achieving the maximum Wiener index is no-zig-zag (i.e. has a core vertex).
\end{conjecture}

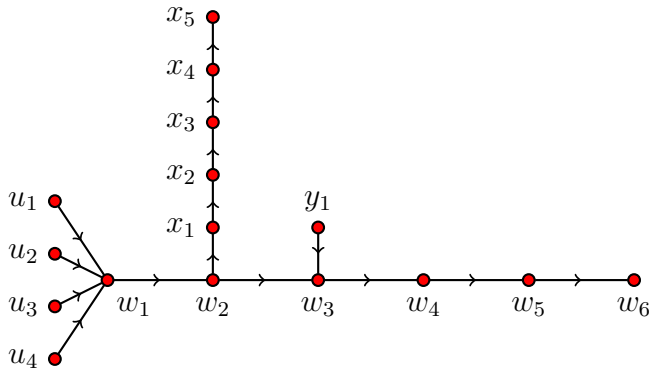
\begin{figure}[ht!]
\begin{center}
\begin{tikzpicture}[scale=1.4,style=thick]


\node [My Style, name=w1]   at (0,0) {};
\node [My Style, name=w2]   at (1,0) {};
\node [My Style, name=w3]   at (2,0) {};
\node [My Style, name=w4]   at (3,0) {};
\node [My Style, name=w5]   at (4,0) {};
\node [My Style, name=w6]   at (5,0) {};

\node [My Style, name=y1]   at (2,0.5) {};

\node [My Style, name=x1]   at (1,0.5) {};
\node [My Style, name=x2]   at (1,1) {};
\node [My Style, name=x3]   at (1,1.5) {};
\node [My Style, name=x4]   at (1,2) {};
\node [My Style, name=x5]   at (1,2.5) {};

\node [My Style, name=u1]   at (-0.5,-0.75) {};
\node [My Style, name=u2]   at (-0.5,-0.25) {};
\node [My Style, name=u3]   at (-0.5,0.25) {};
\node [My Style, name=u4]   at (-0.5,0.75) {};

\draw[middlearrow={>}] (u1) -- (w1);
\draw[middlearrow={>}] (u2) -- (w1);
\draw[middlearrow={>}] (u3) -- (w1);
\draw[middlearrow={>}] (u4) -- (w1);

\draw[middlearrow={>}] (w1) -> (w2);
\draw[middlearrow={>}] (w2) -- (w3); 
\draw[middlearrow={>}] (w3) -- (w4);
\draw[middlearrow={>}] (w4) -- (w5);
\draw[middlearrow={>}] (w5) -- (w6);

\draw[middlearrow={>}] (x1) -> (x2);
\draw[middlearrow={>}] (x2) -- (x3); 
\draw[middlearrow={>}] (x3) -- (x4);
\draw[middlearrow={>}] (x4) -- (x5);

\draw[middlearrow={>}] (w2) -- (x1);

\draw[middlearrow={>}] (y1) -- (w3);

\draw (-0.8,-0.75) node {${u_4}$};
\draw (-0.8,-0.25) node {${u_3}$};
\draw (-0.8,0.25) node {${u_2}$};
\draw (-0.8,0.75) node {${u_1}$};

\draw (0.25,-0.25) node {${w_1}$};
\draw (1,-0.25) node {${w_2}$};
\draw (2,-0.25) node {${w_3}$};
\draw (3,-0.25) node {${w_4}$};
\draw (4,-0.25) node {${w_5}$};
\draw (5,-0.25) node {${w_6}$};

\draw (2,0.75) node {${y_1}$};

\draw (0.7,0.5) node {${x_1}$};
\draw (0.7,1) node {${x_2}$};
\draw (0.7,1.5) node {${x_3}$};
\draw (0.7,2) node {${x_4}$};
\draw (0.7,2.5) node {${x_5}$};

\end{tikzpicture}
\end{center}
\caption{A no-zig-zag tree $T_6$.}

\label{peter}
\end{figure}

It was supported by showing that it holds for trees on at most $10$ vertices, subdivision of stars, and trees constructed from two stars whose central vertices are connected by a path. Furthermore, since it is reasonable to expect that an orientation of a tree
maximizing the Wiener index also maximizes the number of pairs of vertices $(u,v)$ between which there exists a path,  Conjecture \ref{conj:zig-zag} is supported also by a result of Henning and Oellermann \cite{henol}. They proved that if $T$ is a tree and $D$ is an orientation  of $T$ that maximizes the number of ordered pairs $(u,v)$ of vertices of $D$ for which there exists a $(u,v)$-path in
$D$, then $D$ contains a core vertex. However, Li and Wu \cite{LiWu} 
constructed a tree of order $85$ contradicting Conjecture \ref{conj:zig-zag}. Independently, Dankelmann \cite{dankelmann22} found an infinite family of counter-examples. For $k\in \mathbb{N}$, where $k$ is a multiple of $3$, let $T_k$ be the tree obtained from a path of order $k$ with vertices $w_1,w_2,\ldots,w_k$, by connecting vertices $u_1,u_2,\ldots, u_{{k}^2/9}$ to $w_1$, connecting $x_1$ from the path $x_1x_2x_3x_4x_5$ to $w_2$, and a single vertex $y_1$ to $w_3$. Now let $D_k$ be the orientation of $T_k$ such the edges of the path $w_1w_2\ldots w_k$ are oriented towards $w_k$, each edge $u_iw_1$ is oriented towards $w_1$, the
edges of the path $x_1x_2x_3x_4x_5$ are oriented towards $x_5$, and the edge $y_1w_3$ is oriented towards $w_3$, see Figure~\ref{peter} for an example. Observe that the edges of the $(x_5,y_1)$-path change their direction
twice as the path is traversed, thus $D_k$ is a zig-zag orientation. Dankelmann proved that if $k$ is sufficiently large, then $D_k$ and its converse (i.e., a digraph obtained by reversing the
direction of every arc in $D_k$) are the only orientations of $T_k$ that maximize the Wiener index, which contradicts Conjecture \ref{conj:zig-zag}.

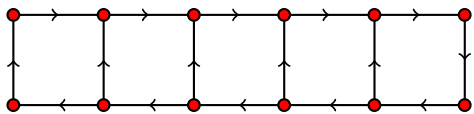
\begin{figure}[ht!]
\begin{center}
\begin{tikzpicture}[scale=1.2,style=thick]


\node [My Style, name=u0]   at (0,0) {};
\node [My Style, name=u1]   at (1,0) {};
\node [My Style, name=u2]   at (2,0) {};
\node [My Style, name=u3]   at (3,0) {};
\node [My Style, name=u4]   at (4,0) {};
\node [My Style, name=u5]   at (5,0) {};

\node [My Style, name=v0]   at (0,1) {};
\node [My Style, name=v1]   at (1,1) {};
\node [My Style, name=v2]   at (2,1) {};
\node [My Style, name=v3]   at (3,1) {};
\node [My Style, name=v4]   at (4,1) {};
\node [My Style, name=v5]   at (5,1) {};

\draw[middlearrow={>}] (u5) -- (u4); \draw[middlearrow={>}] (u4) -> (u3); \draw[middlearrow={>}] (u3) -> (u2); \draw[middlearrow={>}] (u2) -> (u1);\draw[middlearrow={>}] (u1) -> (u0);

\draw[middlearrow={>}] (v0) -- (v1); \draw[middlearrow={>}] (v1) -> (v2); \draw[middlearrow={>}] (v2) -> (v3); \draw[middlearrow={>}] (v3) -> (v4);\draw[middlearrow={>}] (v4) -> (v5);

\draw[middlearrow={>}] (u0) -- (v0); 

\draw[middlearrow={>}] (u1) -- (v1);
\draw[middlearrow={>}] (u2) -- (v2);
\draw[middlearrow={>}] (u3) -- (v3);
\draw[middlearrow={>}] (u4) -- (v4);
\draw[middlearrow={>}] (v5) -- (u5);

\end{tikzpicture}
\end{center}
\caption{An orientation of the ladder $P_6\Box P_2$ with the maximum Wiener index.}

\label{ladder}
\end{figure}

The Cartesian product $P_m\Box P_n$ of paths on $m$ and $n$ vertices, respectively, is called the \textit{grid} and is denoted by $G_{m,n}$. If $m=2$, it is a called the ladder graph $L_n$.
Kraner \v{S}umenjak  et al.~\cite{ladder21} proved a conjecture from \cite{KST-dir-survey} by showing  that the maximum Wiener index of a digraph whose underlying graph is $L_n$ is $(8n^3+3n^2-5n+6)/3$, and is obtained for the orientation presented in Figure \ref{ladder}. In addition, they proved a lower bound for $W_{\max}(G\Box H)$ for general graphs $G$ and $H$, and posed a question regarding its sharpness.
Let $\tau(G)=\sum_{x\in V(G)} \sigma(x)$, where $\sigma(x)$ denotes the number of vertices $x' \in V(G)$ for which there is a path from $x$ to $x'$ in $G$.

\begin{theorem}\label{prod}
For any  graphs $G$ and $H$, 
$$W_{\rm max}(G\Box H) \geq  W_{\rm max}(G) \tau(H)+W_{\rm max}(H)|V(G)|^2\,.$$
\end{theorem}
 
\begin{problem}
Is the bound given in Theorem~\ref{prod} sharp?  Find a sharp lower  bound.
\end{problem}

Another problem from \cite{ladder21} concerns a comparison of the maximum Wiener index of an orientation of $G$ with the Wiener index  of the undirected graph $G$. 
 
\begin{problem}
Find functions $f$ and $g$ so that  $f(W(G))\leq W_{\rm max}(G)\leq g(W(G))$ for all graphs $G$. In particular, can $f$ and $g$
be linear functions?
\end{problem}

\begin{figure}[h!]
\begin{center}
\begin{tikzpicture}[scale=1.2,style=thick]

\node [My Style, name=11]   at (1,1) {};
\node [My Style, name=12]   at (1,2) {};
\node [My Style, name=13]   at (1,3) {};
\node [My Style, name=14]   at (1,4) {};

\node [My Style, name=21]   at (2,1) {};
\node [My Style, name=22]   at (2,2) {};
\node [My Style, name=23]   at (2,3) {};
\node [My Style, name=24]   at (2,4) {};

\node [My Style, name=31]   at (3,1) {};
\node [My Style, name=32]   at (3,2) {};
\node [My Style, name=33]   at (3,3) {};
\node [My Style, name=34]   at (3,4) {};

\node [My Style, name=41]   at (4,1) {};
\node [My Style, name=42]   at (4,2) {};
\node [My Style, name=43]   at (4,3) {};
\node [My Style, name=44]   at (4,4) {};

\node [My Style, name=51]   at (5,1) {};
\node [My Style, name=52]   at (5,2) {};
\node [My Style, name=53]   at (5,3) {};
\node [My Style, name=54]   at (5,4) {};

\node [My Style, name=61]   at (6,1) {};
\node [My Style, name=62]   at (6,2) {};
\node [My Style, name=63]   at (6,3) {};
\node [My Style, name=64]   at (6,4) {};

\draw[middlearrow={>}] (11) -- (12); 
\draw[middlearrow={>}] (12) -- (13);  
\draw[middlearrow={>}] (13) -- (14);

\draw[middlearrow={>}] (21) -- (22); 
\draw[middlearrow={>}] (22) -- (23);  
\draw[middlearrow={>}] (23) -- (24);

\draw[middlearrow={>}] (31) -- (32); 
\draw[middlearrow={>}] (32) -- (33);  
\draw[middlearrow={>}] (33) -- (34); 

\draw[middlearrow={>}] (41) -- (42); 
\draw[middlearrow={>}] (42) -- (43);  
\draw[middlearrow={>}] (43) -- (44); 

\draw[middlearrow={>}] (51) -- (52); 
\draw[middlearrow={>}] (52) -- (53);  
\draw[middlearrow={>}] (53) -- (54);

\draw[middlearrow={>},blue,line width=0.5mm] (62) -- (61); 
\draw[middlearrow={>},blue,line width=0.5mm] (63) -- (62);  
\draw[middlearrow={>},blue,line width=0.5mm] (64) -- (63);

\draw[middlearrow={>}] (61) -- (51); 
\draw[middlearrow={>}] (51) -> (41);
\draw[middlearrow={>}] (41) -- (31); 
\draw[middlearrow={>}] (31) -- (21);  
\draw[middlearrow={>}] (21) -- (11);

\draw[middlearrow={>}] (62) -- (52); 
\draw[middlearrow={>}] (52) -> (42);
\draw[middlearrow={>}] (42) -- (32); 
\draw[middlearrow={>}] (32) -- (22);  
\draw[middlearrow={>}] (22) -- (12);

\draw[middlearrow={>}] (63) -- (53); 
\draw[middlearrow={>}] (53) -> (43);
\draw[middlearrow={>}] (43) -- (33); 
\draw[middlearrow={>}] (33) -- (23);  
\draw[middlearrow={>}] (23) -- (13); 

\draw[middlearrow={>},blue,line width=0.5mm] (54) -- (64); 
\draw[middlearrow={>},blue,line width=0.5mm] (44) -> (54);
\draw[middlearrow={>},blue,line width=0.5mm] (34) -- (44); 
\draw[middlearrow={>},blue,line width=0.5mm] (24) -- (34);  
\draw[middlearrow={>},blue,line width=0.5mm] (14) -- (24);

\draw (1,0.5) node {${(4,1)}$};
\draw (2,0.5) node {${(4,2)}$};
\draw (6,0.5) node {${(4,6)}$};
\draw (1,4.5) node {${(1,1)}$};
\draw (2,4.5) node {${(1,2)}$};
\draw (6,4.5) node {${(1,6)}$};


\node [My Style, name=11a]   at (8,1) {};
\node [My Style, name=12a]   at (8,2) {};
\node [My Style, name=13a]   at (8,3) {};
\node [My Style, name=14a]   at (8,4) {};

\node [My Style, name=21a]   at (9,1) {};
\node [My Style, name=22a]   at (9,2) {};
\node [My Style, name=23a]   at (9,3) {};
\node [My Style, name=24a]   at (9,4) {};

\node [My Style, name=31a]   at (10,1) {};
\node [My Style, name=32a]   at (10,2) {};
\node [My Style, name=33a]   at (10,3) {};
\node [My Style, name=34a]   at (10,4) {};

\node [My Style, name=41a]   at (11,1) {};
\node [My Style, name=42a]   at (11,2) {};
\node [My Style, name=43a]   at (11,3) {};
\node [My Style, name=44a]   at (11,4) {};

\node [My Style, name=51a]   at (12,1) {};
\node [My Style, name=52a]   at (12,2) {};
\node [My Style, name=53a]   at (12,3) {};
\node [My Style, name=54a]   at (12,4) {};

\node [My Style, name=61a]   at (13,1) {};
\node [My Style, name=62a]   at (13,2) {};
\node [My Style, name=63a]   at (13,3) {};
\node [My Style, name=64a]   at (13,4) {};

\draw[middlearrow={>},blue,line width=0.5mm] (11a) -- (12a); 
\draw[middlearrow={>},blue,line width=0.5mm] (12a) -- (13a);  
\draw[middlearrow={>},blue,line width=0.5mm] (13a) -- (14a);

\draw[middlearrow={>},blue,line width=0.5mm] (22a) -- (21a); 
\draw[middlearrow={>},blue,line width=0.5mm] (23a) -- (22a);  
\draw[middlearrow={>}] (23a) -- (24a);

\draw[middlearrow={>},blue,line width=0.5mm] (31a) -- (32a); 
\draw[middlearrow={>},blue,line width=0.5mm] (32a) -- (33a);  
\draw[middlearrow={>}] (33a) -- (34a); 

\draw[middlearrow={>},blue,line width=0.5mm] (42a) -- (41a); 
\draw[middlearrow={>},blue,line width=0.5mm] (43a) -- (42a);  
\draw[middlearrow={>}] (43a) -- (44a); 

\draw[middlearrow={>},blue,line width=0.5mm] (51a) -- (52a); 
\draw[middlearrow={>},blue,line width=0.5mm] (52a) -- (53a);  
\draw[middlearrow={>}] (53a) -- (54a);

\draw[middlearrow={>},blue,line width=0.5mm] (62a) -- (61a); 
\draw[middlearrow={>},blue,line width=0.5mm] (63a) -- (62a);  
\draw[middlearrow={>},blue,line width=0.5mm] (64a) -- (63a);

\draw[middlearrow={>},blue,line width=0.5mm] (61a) -- (51a); 
\draw[middlearrow={>}] (41a) -> (51a); %
\draw[middlearrow={>},blue,line width=0.5mm] (41a) -- (31a); 
\draw[middlearrow={>}] (21a) -- (31a);  %
\draw[middlearrow={>},blue,line width=0.5mm] (21a) -- (11a);

\draw[middlearrow={>}] (52a) -- (62a); %
\draw[middlearrow={>}] (42a) -> (52a); %
\draw[middlearrow={>}] (32a) -- (42a); %
\draw[middlearrow={>}] (22a) -- (32a);  %
\draw[middlearrow={>}] (12a) -- (22a); %

\draw[middlearrow={>}] (53a) -- (63a); 
\draw[middlearrow={>},blue,line width=0.5mm] (53a) -> (43a);
\draw[middlearrow={>}] (33a) -- (43a); 
\draw[middlearrow={>},blue,line width=0.5mm] (33a) -- (23a);  
\draw[middlearrow={>}] (13a) -- (23a); 

\draw[middlearrow={>},blue,line width=0.5mm] (54a) -- (64a); 
\draw[middlearrow={>},blue,line width=0.5mm] (44a) -> (54a);
\draw[middlearrow={>},blue,line width=0.5mm] (34a) -- (44a); 
\draw[middlearrow={>},blue,line width=0.5mm] (24a) -- (34a);  
\draw[middlearrow={>},blue,line width=0.5mm] (14a) -- (24a);

\draw (8,0.5) node {${(4,1)}$};
\draw (9,0.5) node {${(4,2)}$};
\draw (13,0.5) node {${(4,6)}$};
\draw (8,4.5) node {${(1,1)}$};
\draw (9,4.5) node {${(1,2)}$};
\draw (13,4.5) node {${(1,6)}$};

\end{tikzpicture}
\end{center}
\caption{Two orientations, $D_{4,6}$ (left) and $C_{4,6}$ (right), of $P_4\Box P_6$.}

\label{mesh}
\end{figure}
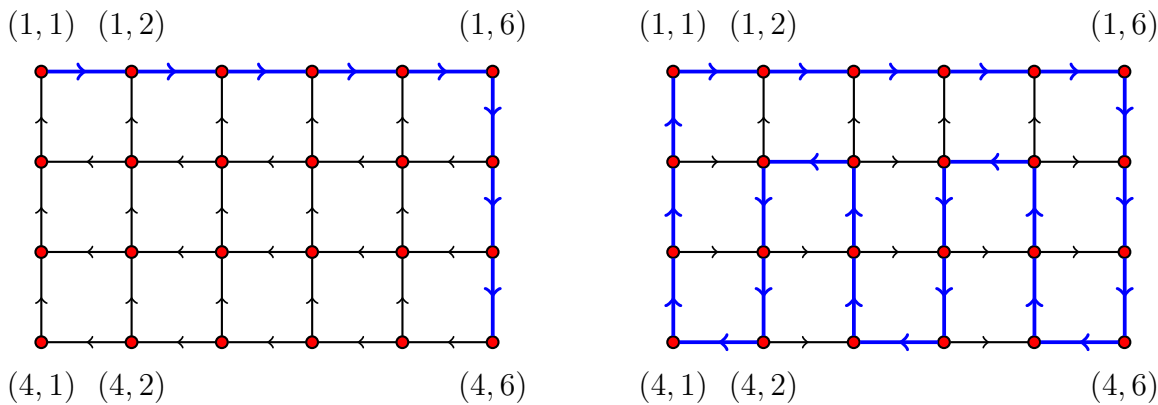

Note that the orientation of $L_n$ in Figure \ref{ladder} is
obtained when all layers isomorphic to one factor are directed paths directed in the same way, except one which is a directed path directed in the opposite way.
Kraner \v{S}umenjak et al.~considered the following natural generalization of this orientation to general grids.
Let ${D}_{m,n}$ be the orientation of $G_{m,n}$ with all
$P_m$-layers oriented up except the last $P_m$-layer which is oriented down,
and all $P_n$-layers oriented to the left except the first $P_n$-layer which
is oriented to the right, see the left graph in Figure \ref{mesh}.

The authors of \cite{ladder21} conjectured that for every $m,n\geq 2$, it holds $W_{\rm max}(G_{m,n})=W({D}_{m,n})$. However, it turns out that a comb-like orientation has significantly bigger Wiener index. Let ${C}_{m,n}$ be an orientation of $G_{m,n}$ in
which the top $P_n$-layer is directed to the right and this layer is
completed to a directed Hamiltonian cycle $C$ in a zig-zag way as shown by blue
arrows on the right graph in Figure~\ref{mesh}.
Moreover, the other edges are directed in such a way that they do not
shorten directed blue path starting at the vertex $(1,1)$.
Of course, ${C}_{m,n}$ exists only if $n$ is even. In \cite{KSgrid}
it was shown that if $n\ge 4$ is even, and $m\ge 3$, then
$W({C}_{m,n})>W({D}_{m,n})$, and further observations led the authors to the following problem.

\begin{problem}
Find the biggest possible constant $c$, such that 
$W_{\max}(G_{m,n})\ge c(mn)^3+o\big((mn)^3\big)$.
\end{problem}

To sum up, the following is still open.

\begin{problem}
Find an orientation of $G_{m,n}$ with the maximum Wiener index.
\end{problem}

The authors think the above problem might be difficult as the extremal graphs in the cases $m=3$ and $n\in \{4,5,6\}$ do not have any
obvious simple property, but they are strongly connected. Thus they ask
the following.

\begin{question}
Let $M_{m,n}$ be an orientation of $G_{m,n}$ with the maximum Wiener index.
Is $M_{m,n}$ strongly connected?
\end{question}

%
%
\section{Maximum Winer index of graphs with prescribed diameter}

Recall that the \textit{eccentricity} of a vertex in a connected graph $G$ is the maximum distance between this vertex and any other vertex of $G$, and the maximum eccentricity is the graph \textit{diameter}. 
Similarly, the {\em radius} of $G$, denoted by ${\rm rad}(G)$,
is the minimum graph eccentricity.
In 1984 Plesn{\'\i}k identified graphs as well as digraphs with a given diameter that minimize the Wiener index (see also \cite{CambieDiam} for a recent alternative proof), and posed the opposite problem regarding the maximum \cite{P}.

\begin{problem}
\label{prob-diam}
What is the maximum Wiener index among graphs of order $n$ and diameter $d$?
\end{problem}

In general this question remains unsolved, but there has been progress and important results were obtained. First, Wang and Guo~\cite{wang2} determined the trees with maximum Wiener index among trees of order $n$ 
and diameter $d$ for some special values of $d$, $2\leq d\leq 4$ or $n-3\leq d \leq n-1$. Mukwembi and Vetr{\'\i}k  \cite{muk} independently considered trees with the diameter up to $6$ and gave asymptotically sharp upper bounds.

DeLaVi{\~n}a and Waller \cite{DeLa} posed a conjecture with additional restrictions in Problem \ref{prob-diam}.

\begin{conjecture}
\label{LaVi}
Let $G$ be a graph with diameter $d>2$
and order $2d + 1$. Then $W(G) \le W(C_{2d+1})$, where $C_{2d+1}$
denotes the cycle of length $2d+1$.
\end{conjecture}

\begin{figure}[h]
	\begin{center}
		\begin{tikzpicture}[scale=0.8,style=thick]
		\node [My Style, name=a]   at (0.2,0.2) {};		
		\node [My Style, name=b]   at (1,1) {};
		\node [My Style, name=c]   at (2.5,2.5) {};		
		\node [My Style, name=d]   at (3.3,3.3) {};
		\node [My Style, name=e]   at (1,0) {};
		\node [My Style, name=f]   at (2.5,0) {};
		\node [My Style, name=g]   at (0,1) {};
		\node [My Style, name=h]   at (2.5,1) {};
		\node [My Style, name=i]   at (3.5,1) {};
		\node [My Style, name=j]   at (0,2.5) {};
		\node [My Style, name=k]   at (1,2.5) {};
		\node [My Style, name=l]   at (3.5,2.5) {};
		\node [My Style, name=m]   at (0.2,3.3) {};
		\node [My Style, name=n]   at (1,3.5) {};
		\node [My Style, name=o]   at (2.5,3.5) {};

		\draw[](a)--(b)--(c)--(d);
		\draw[](m)--(k)--(h);
		\draw[](j)--(k)--(n);
		\draw[](o)--(c)--(l);
		\draw[](g)--(b)--(e);
		\draw[](f)--(h)--(i);
		\draw[](k)--(c)--(h)--(b)--(k);
	
\end{tikzpicture}
\end{center}
\caption{The graph $K_4^{15}$.}
\label{knc}
\end{figure}
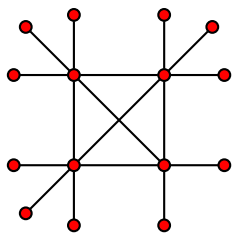

Sun et al.~\cite{sun19} considered general small-diameter and large-diameter graphs. They observed that if $G$ is a graph on $n$ vertices with diameter equal to $2$, then the maximum Wiener index is attained by the star $S_n$. For diameter $3$ they proposed a conjecture, that the extremal graph is isomorphic to $K_n^c$, which is a graph of order $n$ that consists of a complete graph on $c$ vertices and has the rest of the vertices attached to these $c$ vertices as uniformly as possible (meaning that each of the $c$ vertices of the complete graph has either $\lfloor (n-c)/c \rfloor$ or 
$\lceil (n-c)/c \rceil$ pendant vertices attached, see Figure \ref{knc} where 
$K_4^{15}$ is depicted.

\begin{conjecture}
\label{conj:d3}
Let $G$ be a graph on $n$ vertices with diameter equal to $3$. Then
$W(G)\leq W(K_n^c)$ where $c=\Bigl\lfloor \,\sqrt{\frac{n^2}{2(n-1)}} \,\,\Bigr\rfloor$ or $c=\Bigl\lceil \, \sqrt{\frac{n^2}{2(n-1)}} \,\,\Bigr\rceil$.
\end{conjecture}


To explain the results pertaining to trees and a conjecure on general graphs with diameter $4$, we need the following definition. Let $k = \lfloor \sqrt{n-1} \rfloor$. For $k^2 + k \geq n - 1$ we denote by $T_n$ the rooted tree on $n$ vertices in which the root has degree $k$, $n-k^2-1$ of its neighbours are of degree $k+1$ and the rest of them of degree $k$. When $k^2 + k \leq n - 1$ let $T_n'$ denote the rooted tree on $n$ vertices in which the root has degree $k+1$, $n-k^2-k-1$ of its neighbours are of degree $k+1$ and the rest of them of degree $k$. Wang and Guo \cite{wang2} gave a complete description of trees with diameter $4$ that maximize the Wiener index.

\begin{theorem}
\label{thm:d4}
    Let $T$ be a tree on $n$ vertices with diameter $4$ and let $k = \lfloor \sqrt{n-1} \rfloor$. Then the following holds:

    \begin{itemize}
        \item if $k^2 + k > n - 1$, then $W(T) \leq W(T_n)$, with equality holding only when $T \cong T_n$;
        \item if $k^2 + k < n - 1$, then $W(T) \leq W(T_n')$, with equality holding only when $T \cong T_n'$;
        \item if $k^2 + k = n - 1$, then $W(T) \leq W(T_n) = W(T_n')$, with equality holding only when $T \cong T_n$ or $T \cong T_n'$.
    \end{itemize}
\end{theorem}

The authors of~\cite{sun19} suspect that the extremal graphs from the theorem above are extremal also for general graphs.

\begin{conjecture}
\label{conj:d4}
The trees $T_n$ and $T_n'$ remain the unique optima in the class of graphs of diameter $4$ on $n$ vertices as it is described in Theorem~\ref{thm:d4} with the only exception of $n = 9$, in which case $C_9$ is also an optimal graph.
\end{conjecture}

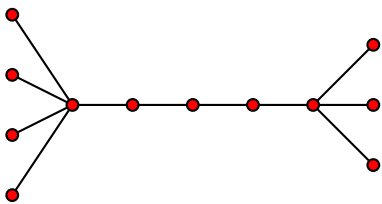
\begin{figure}[h]
	\begin{center}
		\begin{tikzpicture}[scale=0.8,style=thick]
		\node [My Style, name=a]   at (0,3) {};		
		\node [My Style, name=b]   at (0,2) {};
		\node [My Style, name=c]   at (0,1) {};		
		\node [My Style, name=d]   at (0,0) {};
		\node [My Style, name=e]   at (1,1.5) {};
		\node [My Style, name=f]   at (2,1.5) {};
		\node [My Style, name=g]   at (3,1.5) {};
		\node [My Style, name=h]   at (4,1.5) {};
		\node [My Style, name=i]   at (5,1.5) {};
		\node [My Style, name=j]   at (6,2.5) {};
		\node [My Style, name=k]   at (6,1.5) {};
		\node [My Style, name=l]   at (6,0.5) {};

		\draw[](k)--(i)--(l);
		\draw[](c)--(e);
		\draw[](b)--(e)--(d);
		\draw[](a)--(e)--(f)--(g)--(h)--(i)--(j);
	
\end{tikzpicture}
\end{center}
\caption{Double broom $D(12,4,3)$.}
\label{broom}
\end{figure}

An interested reader is referred to \cite{sun19} for computer results supporting Conjectures \ref{LaVi}, \ref{conj:d3} and \ref{conj:d4}. The role of extremal graphs in the case of large-diameter graphs play the so called \textit{double brooms}, i.e.~graphs consisting of a path on $n-a-b$ vertices together with $a$  leaves adjacent to one of its endvertices and $b$ leaves adjacent to the other endvertex (see Figure \ref{broom} for an example).

\begin{theorem}
\label{thm:n-c}
Let $G$ be a graph of order $n$ and diameter $n-c$, where $c \geq 1$ is a constant and $n$ is large enough relative to $c$. Then $W(G) \leq W(D(n, \lfloor (c+1)/2 \rfloor, \lceil (c+1)/2 \rceil))$ with equality if and only if $G \cong D(n, \lfloor (c+1)/2 \rfloor, \lceil (c+1)/2 \rceil)$.
\end{theorem}

Further details on diameters $n-3$ and $n-4$ can be found in \cite{sun19}.
A different approach to Problem \ref{prob-diam} was recently used by Cambie \cite{CambieDiam} who gave asymptotically sharp upper bounds for Wiener index.
As the main first step towards the proof of his result he constructed an almost extremal graph, in which there are many pairs of vertices which are of distance $d$ from each other. This is achieved by having many subtrees with many leaves, and, when the diameter is even, combining them into one tree. When the diameter is odd, a central clique is used so that the distance between leaves of different subtrees are of distance $d$.
Now if we take two vertices at random, the probability that both vertices are leaves is large since the number of leaves is large. Similarly, since we have many subtrees, the probability that both leaves are in different subtrees is large. Hence the probability that two vertices are at maximal distance is large, implying that the average distance is close to $d$. The above is a foundation of the following asymptotic 
solution to the problem of Plesn\'{i}k.

\begin{theorem}
There exist positive constants $c_1$ and $c_2$ such that for any $d\geq 3$ the following holds. The maximum Wiener index among all graphs of diameter $d$ and order $n$ is between $d-c_1\frac{d^{3/2}}{\sqrt{n}}$ and $d-c_2\frac{d^{3/2}}{\sqrt{n}}$, i.e.~it is of the form $d-\Theta\left( \frac{d^{3/2}}{\sqrt{n}} \right)$.
\end{theorem}

In addition, Cambie \cite{CambieDiam} gives slightly stronger upper bound for trees, by which he extends a result of Mukwembi and Vetr\'{i}k \cite{muk}. Moreover, the results he obtained lead him to the following question.

\begin{question}
For even $d$ and large $n$, are the graphs of order $n$ and diameter $d$ with the largest Wiener index all trees?
\end{question}

Digraphs were considered in \cite{CambieDiam} as well, where the problem of Plesn\'{i}k is solved exactly if the order is large comparing to the diameter. For the sake of completeness we also mention that trees of order $n$ and diameter $d$ with the minimum Wiener index were presented in \cite{liupan}.

\medskip

Having in mind the close relationship between the diameter
and the radius of a connected graph, ${\rm rad}(G)\leq {\rm diam}(G) \leq 2\, {\rm rad}(G)$, it is natural to consider
the above problems with radius instead of diameter. Chen et al.~\cite{Chen} posed the following question.

\begin{problem} \label{n-r} What is the maximum Wiener index among
graphs of order $n$ and radius $r$?
\end{problem}

They succeeded to characterize graphs with the maximum Wiener
index among all graphs of order $n$ with radius $2$. 
Das and Nadjafi-Arani \cite{das17} gave an upper bound on Wiener index
of trees and graphs in terms of number of vertices $n$ and radius $r$. In addition, they presented an upper bound on the
Wiener index in terms of order, radius and maximum degree
of trees and of graphs. 
The authors concluded that these results are not enough to solve Problem \ref{n-r}.  Stevanovi\'{c} et al.~\cite{stevan20} provide examples obtained by computer experiments, which suggest that a simple characterization of the structure of trees with maximum Wiener index among trees with a given number of vertices and radius will probably be out of our reach in some foreseeable future.

Analogous problem for the minimum Wiener index was posed by You and Liu \cite{You}.

\begin{problem} What is the minimum Wiener index among all graphs of order $n$ and radius $r$?
\end{problem}

If $r\in \{1,2\}$, the extremal graphs attaining the minimal total distance among all graphs of order $n$ are easily characterized: they are complete graphs when $r = 1$, complete graphs minus a
maximum matching when $r = 2$ and $n$ is even, and complete graphs minus a maximum matching and an additional edge adjacent to the vertex not in the maximum matching, when $r = 2$ and $n$ is odd. 

A conjecture for $n\geq 3$ was posed by Chen et al. \cite{Chen}.
The notation $G_{n,r,s}$, where $n,r$ and $s$ are positive integers such that $n\geq 2r$, $r\ge 3$, and $n-2r+1\ge s\ge 1$, stands for the graph obtained in the following way: let $v_1, v_2, v_3$ and $v_4$ be four
consecutive vertices on a $2r$-cycle.  Replace $v_2$ with a clique of order
$s$, replace $v_3$ with a clique of order $n-2r+2-s$, join each
vertex of one clique to all vertices of the other clique, join $v_1$ to
the all vertices of $K_s$, and join $v_4$ to all vertices of $K_{n-2r+2-s}$. Notice that the resulting graph has $n$ vertices and radius $r$, and 
$W(G_{n,r,s}) = W(G_{n,r,s'})$ for any $s,s'\in \{1,\ldots,r-1\}$.

\begin{conjecture}
Let $n$ and $r$ be two positive integers with $n\geq 2r$ and $r\ge 3$. 
For any graph $G$ of order $n$ with radius $r$, $W(G)\geq  W(G_{n,r,1})$. Equality is attained if and only if $G=G_{n,r,s}$ for $s\in \{1,\ldots,r-1\}$.
\end{conjecture}

Cambie showed that the hypercube $Q_3$ is a counterexample to the above conjecture, so it does not hold when $n$ is small, but he demonstrated that the conjecture is true asymptotically, i.e. if the order is sufficiently large compared to the radius~\cite{Cambie21rad}.

\begin{theorem}
For any $r\ge 3$, there exists a value $n_1(r)$ such that for all $n\geq n_1(r)$ it holds that any graph  $G$ of order $n$ with radius $r$ satisfies $W(G)\geq  W(G_{n,r,1})$. Equality holds if and only if $G=G_{n,r,s}$ for $s\in \{1,\ldots,r-1\}$.
\end{theorem}

We refer to~\cite{Cambie21rad} for an analog of this result for directed graphs, and to~\cite{P} for a characterization of digraphs of given order and diameter with the minimum Wiener index.




%
%
\section{\v{S}olt\'{e}s problem and its relaxed variations}

An interesting question regarding the Wiener index is to study how Wiener index is affected by small changes in a graph. Clearly, by removing an edge Wiener index is increased. On the other hand, the effect of deleting a vertex is far from obvious, and it was first studied by \v{S}olt\'{e}s. In his paper from 1991, \v{S}olt\'{e}s posed the following problem~\cite{soltes}.

\begin{problem}\label{solt}
Find all graphs $G$ in which the equality $W(G) = W(G-v)$ holds for all $v\in V(G)$.
\end{problem}

Therefore, if for a vertex $v$ in a graph $G$ it holds that $W(G) = W(G-v)$, we say that $v$ satisfies the \textit{\v{S}olt\'{e}s property} in $G$, and a graph in which every vertex satisfies the \v{S}olt\'{e}s property is referred to as a \textit{\v{S}olt\'{e}s graph}. The only known \v{S}olt\'{e}s graph so far is the cycle on $11$ vertices. The above problem appears to be difficult, thus in subsequent studies relaxed variations were considered.
The authors of \cite{Majst-1}
showed that the class of graphs for which the Wiener index does not change when a particular vertex is removed is rich, even when restricted to unicyclic graphs with  fixed length of the cycle. More precisely:
\begin{itemize}
\item there is a unicyclic graph $G$ on $n$ vertices containing a vertex $v$  with $W(G) = W(G-v)$ if and only if $n \geq 9$;
\item there is a unicyclic graph $G$ with a cycle of length $c$ and a vertex satisfying the \v{S}olt\'{e}s property
if and only if $c \geq 5$;
\item for every graph $G$ there are infinitely many graphs $H$ such that $G$ is an induced subgraph of $H$ and $W(H) = W(H-v)$ 
for some $v \in V(H) \setminus V(G)$.
\end{itemize}
If a vertex $v$ has degree $1$ in $G$, then clearly $W(G) > W(G-v)$. In the construction of the above mentioned infinite class of graphs $G$ with a vertex $v$ satisfying the \v{S}olt\'{e}s property the vertex $v$ is of degree $2$. In \cite{Majst-2} the authors extended their research to graphs in which $v$ is of arbitrary degree. They showed that for a fixed positive integer $k\geq 2$ there exist infinitely many graphs $G$ with a vertex $v$ such that $\deg_G(v) = k$ and 
$W(G) = W(G-v)$. Moreover, if $n\geq 7$, there exists an $n$-vertex graph $G$ with a vertex $v$ so that $\deg_G(v) = n-2$  or $\deg_G(v) = n-1$, respectively, and $W(G) = W(G-v)$. By proving the next theorem they showed that dense graphs cannot be a solution of Problem \ref{solt}.

\begin{theorem}
\label{tm:min_d}
If $G$ is an $n$-vertex graph for which $\delta(G)\geq n/2$,
then $W(G)\neq  W(G-v)$ for every $v\in V(G)$.
\end{theorem}

In the results above, removal of one vertex only was considered. So  
the authors proposed the study of graphs $G$ in which a given number of vertices satisfying the \v{S}olt\'{e}s property exist \cite{Majst-2,Majst-3}.

\begin{problem}
\label{infk}
For a given $k$, find (infinitely many) graphs $G$ for which $W(G)=W(G-v_1)=W(G-v_2)=\cdots =W(G-v_k)$ for some distinct vertices $v_1,\ldots,v_k$ in $G$.
\end{problem}

This problem was considered by Bok et al.~\cite{Bok1,Bok2} who
showed the existence of:
\begin{itemize}
\item infinitely many \textit{cactus graphs} (i.e.~graphs in which every edge belongs to at
most one cycle) with exactly $k$ cycles of length at least $7$ that contain exactly $2k$ vertices satisfying the \v{S}olt\'{e}s property; and 
\item infinitely many cactus graphs with exactly $k$ cycles of length $c \in \{5,6\}$ that contain exactly $k$ vertices satisfying the \v{S}olt\'{e}s property. 
\end{itemize}
In addition, they proved that $G$ contains no vertex with the \v{S}olt\'{e}s property if the length of the longest cycle in $G$ is at most $4$.
Another infinite family of graphs satisfying the condition from Problem \ref{infk} was constructed by Hu et al.~\cite{Hu21}. Furthermore, Hu et al.~settled another problem from \cite{Majst-2, Majst-3} by proving that for any $k\geq 2$, there exist infinitely many graphs $G$ such
that $W(G) = W(G-\{v_1,v_2,\ldots,v_k\})$ for some distinct vertices $v_1,v_2,\ldots,v_k\in V(G)$.




	

\medskip

Akhmejanova et.~al \cite{akh} considered a relaxation of the original \v{S}olt\'{e}s problem from another point of view. They asked for graphs with a large proportion of vertices satisfying the \v{S}olt\'{e}s property. More precisely, they defined the function $\Delta_v(G)=W(G)-W(G-v)$. Then $$\frac{|\{v\in V(G);\Delta_v(G)=0\}|}{|V(G)|}$$
is the proportion of vertices satisfying the \v{S}olt\'{e}s property. So Akhmejanova et.~al asked the following. 

\begin{problem}
For a fixed $\alpha \in (0, 1]$ construct an infinite series $S$ of graphs such that for all $G = (V(G),E(G))$
from $S$ the following holds:
$$\frac{|\{v\in V(G);\Delta_v(G)=0\}|}{|V(G)|}\geq \alpha.$$
\end{problem}

\noindent Note that a solution to this problem for $\alpha = 1$ would give an infinite series of solutions to Problem~\ref{solt}. The authors noted that a slight modification of a construction from \cite{Bok1} yields an infinite series of graphs with the proportion of vertices satisfying the \v{S}olt\'{e}s property tending to $\frac{1}{3}$, and improved this constant by finding another two constructions. The first construction contains many $11$-cycles as induced subgraphs: given $k\in \mathbb{N}, k>1$, they defined a graph $B(k)$ on $5k+6$ vertices by taking two vertices and connecting them with $k$ distinct paths of length $6$ and one path of length $5$. It turns out that  for $B(k)$ the proportion of vertices satisfying the \v{S}olt\'{e}s property equals $\frac{2k}{5k+6}$, thus this proportion tends to $\frac{2}{5}$ as $k$ tends to infinity.
Another construction of so called lily-shaped graphs involves graphs that are not $2$-connected and whose proportion tends to $\frac{1}{2}$, see \cite{akh} for details. Furthermore, the authors found a graph with the proportion $\frac{2}{3}$ and expect that there exist an infinite series of graphs with a proportion $\alpha > \frac{1}{2}$, or perhaps even $\alpha$ tending to $1$. Furthermore, they propose the following problems.

\begin{problem}
For a fixed $z \in \mathbb{Z}$, find all graphs $G$, for which the equality $W(G)-W(G-v) = z$ holds for all vertices $v$.
\end{problem}

\begin{problem}
For a fixed $z \in \mathbb{Z}$ and $\alpha \in (0, 1]$, construct an infinite series $S$ of graphs such that for all
$G = (V(G),E(G))$
from $S$ the following inequality takes place:
$$\frac{|\{v\in V(G);\Delta_v(G)=z\}|}{|V(G)|}\geq \alpha.$$
\end{problem}

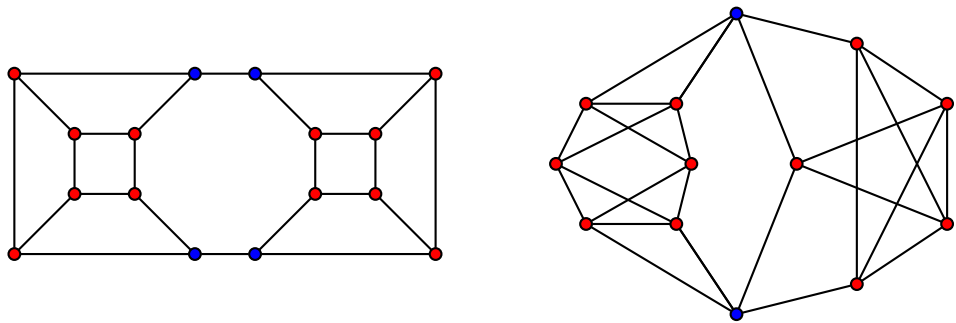
\begin{figure}[h!]
	\begin{center}
		\begin{tikzpicture}[scale=0.4,style=thick]
	
		\node [My Style, name=a1]   at (-18,2) {};		
		\node [My Style, name=b1, fill=blue]   at (-12,2) {};
		\node [My Style, name=c1, fill=blue]   at (-10,2) {};		
		\node [My Style, name=d1]   at (-4,2) {};
		\node [My Style, name=e1]   at (-4,8) {};
		\node [My Style, name=f1, fill=blue]   at (-10,8) {};
		\node [My Style, name=g1, fill=blue]   at (-12,8) {};
		\node [My Style, name=h1]   at (-18,8) {};
		\node [My Style, name=i1]   at (-16,6) {};
		\node [My Style, name=j1]   at (-14,6) {};
		\node [My Style, name=k1]   at (-14,4) {};
		\node [My Style, name=l1]   at (-16,4) {};
		\node [My Style, name=m1]   at (-8,6) {};
		\node [My Style, name=n1]   at (-6,6) {};
		\node [My Style, name=o1]   at (-6,4) {};
		\node [My Style, name=p1]   at (-8,4) {};

		\draw[](a1)--(b1)--(c1)--(d1)--(e1)--(f1)--(g1)--(h1)--(a1);
		\draw[](i1)--(j1)--(k1)--(l1)--(a1);
		\draw[](m1)--(p1)--(o1)--(d1);
		\draw[](f1)--(m1)--(n1)--(e1);
		\draw[](h1)--(i1)--(l1);
		\draw[](j1)--(g1);
		\draw[](k1)--(b1);
		\draw[](n1)--(o1);
		\draw[](p1)--(c1);

		\node [My Style, name=a]   at (0,5) {};		
		\node [My Style, name=b]   at (1,7) {};
		\node [My Style, name=c, fill=blue]   at (6,10) {};		
		\node [My Style, name=d]   at (10,9) {};
		\node [My Style, name=e]   at (13,7) {};
		\node [My Style, name=f]   at (13,3) {};
		\node [My Style, name=g]   at (10,1) {};
		\node [My Style, name=h, fill=blue]   at (6,0) {};
		\node [My Style, name=i]   at (1,3) {};
		\node [My Style, name=j]   at (4.5,5) {};
		\node [My Style, name=k]   at (4,7) {};
		\node [My Style, name=l]   at (4,3) {};
		\node [My Style, name=m]   at (8,5) {};

		\draw[](e)--(m)--(f);
		\draw[](c)--(m)--(h);
		\draw[](g)--(d)--(f);
		\draw[](e)--(g);
		\draw[](a)--(b)--(c)--(d)--(e)--(f)--(g)--(h)--(i)--(a);
		\draw[](c)--(k)--(j)--(l)--(h);
		\draw[](c)--(k)--(a)--(l)--(h);
		\draw[](k)--(b)--(j)--(i)--(l);
	
\end{tikzpicture}
\end{center}
\caption{Regular graphs with blue vertices satisfying the \v{S}olt\'{e}s property.}
\label{solt2}
\end{figure}

In \cite{Majst-2, Majst-3} the problem of finding $k$-regular connected graphs\, $G$ other than\, $C_{11}$ for which the equality $W(G)=W(G-v)$ holds for at least one vertex $v\in V(G)$ was posed. The answer is affirmative, see Figure \ref{solt2} for $3$-regular and $4$-regular graphs with $4$ and $2$, respectively, (blue) vertices satisfying the \v{S}olt\'{e}s property. 
Using computer software and counting cubic graphs of orders $n\leq 26$, Ba\v{s}i\'{c} et al.~\cite{basic} found that cubic graphs of order $12$ or less do not contain \v{S}olt\'{e}s vertices. Cubic graphs with two \v{S}olt\'{e}s vertices first appear at the order $14$ (there are three such graphs), and examples with three and four \v{S}olt\'{e}s vertices appear at the order $16$. Moreover, they proved the following.

\begin{theorem}
There exist infinitely many cubic $2$-connected graphs which contain two
\v{S}olt\'{e}s vertices. 
\end{theorem}

\medskip

In the same paper, graphs where the
ratio between the number of \v{S}olt\'{e}s vertices and the order of the graph is at least $\alpha$ are called $\alpha$-\v{S}olt\'{e}s graphs. So Problem \ref{solt} asks to find all $1$-\v{S}olt\'{e}s graphs.
The authors believe the solution to this problem should be graphs having all vertices of the same degree.

\begin{conjecture}
If $G$ is a \v{S}olt\'{e}s graph, then it is regular.
\end{conjecture}

For a general regular graph $G$, the values $W(G-u)$ and $W(G-v)$ might be significantly different for two different vertices $u$ and $v$ from $G$. It may happen that removal of one vertex
increases the Wiener index, while removal of the other vertex descreases it. However, $W(G-u)$ and $W(G-v)$ are equal if vertices $u$ and $v$ belong to the same vertex
orbit. This led the authors to believe the following.

\begin{conjecture}
If $G$ is a \v{S}olt\'{e}s graph, then $G$ is vertex-transitive.
\end{conjecture}

Further, the authors report that a computer search on publicly available
collections of vertex-transitive graphs did not reveal any $1$-\v{S}olt\'{e}s graphs. All examples of $\frac{1}{3}$-\v{S}olt\'{e}s graphs
are obtained by truncating certain cubic vertex-transitive graphs, and there are no \v{S}olt\'{e}s graphs among vertex-transitive
graphs with less than $48$ vertices. Therefore it is reasonable to consider the following conjectures and a problem.

\begin{conjecture}
If $G$ is a \v{S}olt\'{e}s graph, then $G$ is a Cayley graph.
\end{conjecture}

\begin{problem}
Find an infinite family of cubic vertex-transitive graphs $\{G_i\}_{i=1}^{\infty}$, such that the truncation of $G_i$  is a $\frac{1}{3}$-\v{S}olt\'{e}s graph for all $i\geq 1$.
\end{problem}

\begin{conjecture}
The cycle on eleven vertices is the only \v{S}olt\'{e}s graph.
\end{conjecture}

\section{Wiener index of signed graphs}

A \textit{signed graph} is a pair $(G,\sigma)$ where $G$ is a graph and $\sigma$ is a function from $E(G)$ to $\{-1,1\}$, called
a \textit{signature function} (also called \textit{signing} in the literature). A path $P$ is a $uv$-\textit{path} if its endvertices are $u$ and $v$. If $P$ is a path in $G$ and $\sigma$ is a signature function of $G$ then the notation $\sigma(P)$ stands for the sum $\sum_{e\in P}\sigma (e)$.
For $u,v \in V(G)$ the \textit{signed distance} $d_{G,\sigma}(u,v)$ equals $\min_P |\sigma(P)|$ where the minimum ranges
over all $uv$-paths $P$. Spiro \cite{spiro21} recently introduced the Wiener index $W_{\sigma}(G)$ of the signed graph $(G,\sigma)$ as
$$W_{\sigma}(G)=\displaystyle \sum_{\{u,v\}\subseteq V(G)}d_{G,\sigma}(u,v).$$
If $\sigma$ is a constant function, then $d_{G,\sigma}(u,v) = d(u,v)$, and therefore $W_{\sigma}(G)= W(G)$.
In particular, if $W(G) = W(G-v)$ for all $v\in V(G)$, then there exists a (constant) signature function $\sigma$ of $G$ such that $W_{\sigma}(G) = W_{\sigma}(G-v)$. In this sense the problem of finding signed graphs $(G,\sigma)$ with
$W_{\sigma}(G) = W_{\sigma}(G-v)$ can be viewed as a relaxation of \v{S}olt\'{e}s problem. Note that in the signed setting, it is possible to have
$W_{\sigma}(G)=0$. Spiro used this fact to provide many examples of signed graphs satisfying $W_{\sigma}(G) = W_{\sigma}(G-v)$ for all
$v\in V(G)$, and even with $W_{\sigma}(G) = W_{\sigma}(G-S)$ for any set $S$ of size less than some value $k$. To present his results, a signature function $\sigma$ of a graph $G$ is called \textit{$k$-canceling} if for any set $S \subseteq V(G)$ of
size less than $k$, we have $W_{\sigma}(G-S) = 0$. A graph $G$ is \textit{$k$-canceling} if there exists
a $k$-canceling signature function $\sigma$ of $G$, and graphs with $W_{\sigma}(G) = 0$ are simply referred to as \textsl{canceling graphs}. For instance, a complete graph $K_n$ is $k$-canceling if $n\geq 2k+4$. Furthermore, he proved the following.

\begin{proposition}
Let $G'$ be a bipartite graph with partite sets $U$ and $V$, where$|U|,|V|\geq k + 2$, and minimum
degree at least $k + 1$. Let $G$ be the graph obtained from $G'$ by adding every edge between two vertices of $U$ and every edge between two vertices of $V$. Then G is $k$-canceling.
\end{proposition}

Another family of examples is obtained from the blowups of odd cycles: if $G$ is a graph on $\{v_1,\ldots, v_t\}$, then the \textit{$\{n_1,\ldots, n_t\}$-blowup} of $G$ is defined to be the $t$-partite graph on sets $V_1,\ldots V_t$ with $|V_i|= n_i$ and with $u\in V_i$ and $w\in V_j$ adjacent if and only if $v_i,v_j$ are adjacent in $G$.

\begin{proposition}
Let $G$ be the $(n_1,\ldots, n_{2t+1})$-blowup of a cycle $C_{2t+1}$ with $t\geq 1$. If $n_i \geq 2k$ for all $i$, then $G$ is $k$-canceling.
\end{proposition}

Furthermore, the following holds.
\begin{theorem}
If $n$ is sufficiently large and $G$ is an $n$-vertex graph with minimum degree at
least $\frac{2n}{3}$, then there exists a signature function $\sigma$ of $G$ such that $W_{\sigma}(G) = W_{\sigma}(G-v) = 0$ for all $v \in V(G)$.
\end{theorem}

For necessary conditions for a graph to be canceling and several interesting open questions we refer to \cite{spiro21}. One of the conjectures pertains to the well known fact that in the class of $n$-vertex trees the star $S_n$ and the path $P_n$ are extremal graphs for the Wiener index. Let $(T,\sigma)$ be a signed $n$-vertex tree and let $+$ be the constant signature function that assigns $+1$ to every edge of $P_n$. Then the fact that $W_{\sigma}(T) \leq W_+(P_n)$ follows from the result for the classical Wiener
index since $W_+(P_n) = W(P_n)$. It remains to prove the lower bound.

\begin{conjecture}
If $(T,\sigma)$ is a signed $n$-vertex tree, then
$$W_{\alpha}(P_n) \leq W_{\sigma}(T),$$
where $\alpha$ is the alternating
signature function which assigns the first edge of the path $+1$, the second $-1$, the third $+1$, and so on.
\end{conjecture}

Another possible direction for future study according to Spiro is the 
\textit{minimum signed Wiener index} $W_*(G)=\min_{\sigma}(G)$, where the minimum ranges over all signature functions $\sigma$ of $G$. Note that this concept is analogous to the minimum digraph Wiener index of all orientations of a graph $G$ presented in Section \ref{4}. Spiro proposed a conjecture in which double stars appear as extremal graphs; a \textit{double star} is a tree $T$ in which there exist vertices $x,y \in V(T)$ such that every edge of $T$ has at least one of the vertices $x,y$ as an end vertex. Note that by this definition a star is also a double star.

\begin{conjecture}
If $T$ is an $n$-vertex tree, then
$$W_*(P_n) \leq W_*(T) \leq \displaystyle \max_{D \in \mathcal{D}}W_*(D),$$
where $\mathcal{D}$ is the set of all $n$-vertex double stars.
\end{conjecture}

The conjecture was verified for $n \leq 9$, and noted that it is false if one considers stars instead of double stars. We refer to \cite{spiro21} for more interesting questions related to the presented topic.

%
%
\section{Variable Wiener index vs.~Variable Szeged index}

For an edge $uv$ in a graph, let $n_v(u)$ denote the number of vertices strictly closer to $u$
than $v$, and analogously, let $n_u(v)$ be the number of vertices strictly
closer to $v$ than $u$. In his original paper \cite{Wiener} Wiener observed that the Wiener index of a tree can be computed as the sum of products $n_v(u)\cdot n_u(v)$ over all edges $uv$ in the tree, but this is not the case in general graphs, owing to the fact that shortest paths are typically not unique. By relaxing the condition that the graph is a tree, the Szeged index of a graph $G$ was defined  in~\cite{Gut94,KhaDesKal95} as
$$ 
\Sz(G)= \sum_{uv \in E(G)} n_v(u)\cdot n_u(v).
$$ 

\noindent
Klav\v{z}ar et al.~\cite{KlaRajGut96} proved that $\Sz(G)\ge W(G)$ for every graph $G$, and in~\cite{DobGut95} all graphs for which
the equality holds were classified.

\begin{theorem}
\label{thm:SzW}
For every graph $G$ we have
$\Sz(G)\ge W(G)$,
and equality holds if and only if every block of $G$ is a complete graph. 
\end{theorem}


The {\em variable Wiener index} (also known as the 
{\em generalized Wiener index}) of a graph $G$ is defined as
\begin{equation*}\label{varW}
\displaystyle W^{\alpha}(G)=\sum_{\{u, v\} \subseteq V(G)} d(u,v)^{\alpha},
\end{equation*}
and the {\em variable Szeged index} of a graph $G$ is
\begin{equation*}\label{varSz}
\displaystyle \Sz^{\alpha}(G)=\sum_{uv\in E(G)} \left (n_v(u)\cdot n_u(v)\right)^{\alpha}.
\end{equation*}
Note that in \cite{GutVukZer04} the quantity $\sum_{uv\in E(T)} \left (n_v(u)\cdot n_u(v)\right )^{\alpha}$ was named as the \textsl{variable Wiener index for trees}, but referring to it as 
the \textsl{variable Szeged index} seems to be more natural. 
By Theorem~{\ref{thm:SzW}}, for trees it holds $W(T)=\Sz(T)$.
Using Karamata's inequality Hri\v{n}\'{a}kov\'{a} et al.~\cite{HKS19} proved the following statement.

\begin{theorem}
\label{thm:tree}
Let $T$ be a tree on $n$ vertices.
Then
\begin{enumerate}
\item
$W^{\alpha}(T)\le \Sz^{\alpha}(T)$ if $\alpha>1$,
\item
$W^{\alpha}(T)\ge \Sz^{\alpha}(T)$ if $0\le\alpha<1$.
\end{enumerate}
Moreover, equalities hold if and only if $n=2$.
\end{theorem}

In the case when $\alpha >1$, they extended this result to the class of bipartite graphs.

\begin{theorem}
\label{thm:bipartite}
Let $G$ be a bipartite graph on $n$ vertices and $\alpha>1$.
Then $W^{\alpha}(G)\le\Sz^{\alpha}(G)$ with equality if and
only if $n=2$.
\end{theorem}

If $G$ is a complete graph, we have $\Sz^{\alpha}(G) =\binom{|V(G)|}{2} = W^{\alpha}(G)$ for every $\alpha$.
Note that $\alpha$ is non-negative in the above results. If $\alpha< 0$ then for non-complete graphs we have the following strict inequality \cite{HKS19}.

\begin{proposition}
Let $G$ be a non-complete graph. Then for every $\alpha<0$ we have $\Sz^{\alpha}(G) < W^{\alpha}(G)$.
\end{proposition}

Based on Theorem~{\ref{thm:tree}} and examples provided in \cite{HKS19}, Hri\v{n}\'{a}kov\'{a} et al.~proposed the following conjecture.

\begin{conjecture}\label{strong} 
For every non-complete graph $G$ there is a constant $\alpha_G\in(0,1]$ such that
\begin{eqnarray*}
\Sz^{\alpha}(G) &>& W^{\alpha}(G),\qquad
\mbox{if $\alpha>\alpha_G$},\\
\Sz^{\alpha}(G) &=& W^{\alpha}(G),\qquad
\mbox{if $\alpha=\alpha_G$},\\
\Sz^{\alpha}(G) &<& W^{\alpha}(G),\qquad
\mbox{if $0\leq\alpha<\alpha_G$.}
\end{eqnarray*}
\end{conjecture}

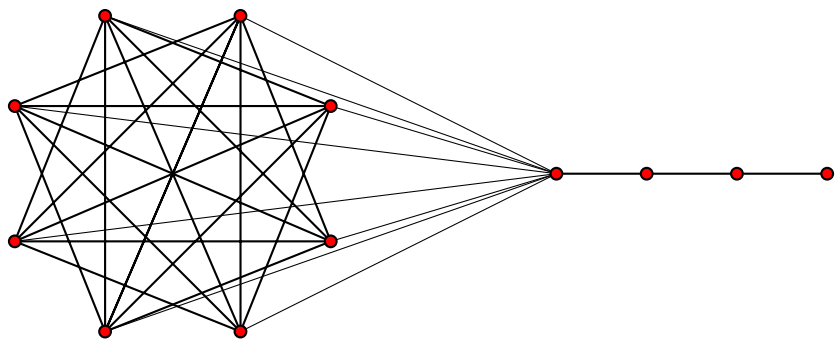
\begin{figure}[ht!]
\begin{center}
\begin{tikzpicture}[scale=1.2,style=thick]

		\node [My Style, name=a]   at (1,2.5) {};		
		\node [My Style, name=b]   at (2.5,2.5) {};
		\node [My Style, name=c]   at (3.5,1.5) {};		
		\node [My Style, name=d]   at (3.5,0) {};
		\node [My Style, name=e]   at (2.5,-1) {};
		\node [My Style, name=f]   at (1,-1) {};
		\node [My Style, name=g]   at (0,0) {};
		\node [My Style, name=h]   at (0,1.5) {};
		\node [My Style, name=i]   at (6,0.75) {};
		\node [My Style, name=j]   at (7,0.75) {};
		\node [My Style, name=k]   at (8,0.75) {};
		\node [My Style, name=l]   at (9,0.75) {};

		\draw[](a)--(c)--(e)--(g)--(a)--(f)--(c)--(h)--(d)--(a)--(e)--(b)--(d)--(g)--(b)--(h)--(f)--(b)--(f);
		\draw[](h)--(e);
		\draw[](c)--(g);
		\draw[](f)--(d);
		\draw[](i)--(j)--(k)--(l);
		
		\draw[thin](a)--(i);
		\draw[thin](b)--(i);
		\draw[thin](c)--(i);
		\draw[thin](d)--(i);
		\draw[thin](e)--(i);
		\draw[thin](f)--(i);
		\draw[thin](g)--(i);
		\draw[thin](h)--(i);

\end{tikzpicture}
\end{center}
\caption{The graph $G_{k,\ell}$ for $k=8$ and $\ell=3$.}

\label{cambie}
\end{figure}

In other words, the conjecture states that for
any non-complete graph there is a critical exponent in $(0,1]$, below which the variable Wiener index is larger and above which the variable Szeged index is larger. As seen above, this holds for trees. 
However, Cambie and Haslegrave \cite{CambHas} found infinitely many counterexamples by constructing a family of graphs $G_{k,\ell}$ as follows: take a complete graph $K_k$, remove a $k$-cycle from it, and connect all its vertices with one endvertex of a path of length $l$, see Figure~\ref{cambie} where $G_{8,3}$ is depicted. 
By fixing a connected non-complete graph $G$, $h(\alpha)=Sz^{\alpha}(G)-W^{\alpha}(G)$ is a continuous function with $h(0)<0$ and $h(1)\geq 0$, which by intermediate value theorem implies that there is at least one value of $\alpha$ for which $h(\alpha)=0$, and at least one such value lies in $(0,1]$. Therefore Conjecture \ref{strong} is equivalent to $\alpha$ being unique, which is not the case for many graphs of the form $G_{k,\ell}$. It turns out that if $k$ is reasonably large, then there exist some corresponding values of $\ell$ having three values of $\alpha$ for which $Sz^{\alpha}(G_{k,\ell})-W^{\alpha}(G_{k,\ell})$ equals $0$.


On the other hand, the authors found further families of graphs for which the statement in Conjecture \ref{strong} does hold. In fact, they showed its validity for almost all graphs.

\begin{theorem}
Conjecture \ref{strong} holds for
\begin{itemize}
\item block graphs,
\item edge-transitive graphs,
\item bipartite graphs,
\item graphs with diameter $2$,
\item graphs with diameter $3$, $n$ vertices and at most $\frac{1}{2}\binom n2$ edges,
\item graphs with $n$ vertices and $m$ edges whenever $m\leq \frac{1}{4}(n^{4/3}-n^{1/3})$.
\end{itemize}
\end{theorem}

They proved also that Conjecture \ref{strong}  holds for almost all random graphs in 2 models of random graphs, see \cite{CambHas}  for more detailed explanation.
Anyway, it is an open problem if there exist graphs $G$, other than
complete ones, for which $|\{\alpha; Sz^{\alpha}(G)-W^{\alpha}(G)=0\}|$ is
larger than 3.
So we have the following problem.

\begin{problem} Let $\mathcal G$ be the class of graphs which contain at
least one block which is not complete.
Is $|\{\alpha; Sz^{\alpha}(G)-W^{\alpha}(G)=0\}|$ bounded for
$G\in\mathcal G$?
If so, what is its maximum value?
\end{problem}

By showing that for every graph $G$, the sequence $(n_v(u)\cdot n_u(v))_{uv\in E(G)}$ majorizes the sequence $(d(u,v))_{u\in V(G)}$,  Cambie and Haslegrave proved that a weaker version of Conjecture \ref{strong} holds. Using a different approach the same result was independently obtained by Kovijani\'{c} Vuki\'{c}evi\'{c} and Bulatovi\'{c} \cite{ZanaBul}.

\begin{theorem}
For every non-complete  graph $G$ and $\alpha> 1$, we have $\Sz^{\alpha}(G)>W^{\alpha}(G)$.
\end{theorem}




\section{Wiener index of apex graphs}

An \textit{apex graph} is a graph that becomes planar by removal of a single vertex. Along these lines a graph $G$ is called an \textit{apex tree} if it contains a vertex $x$ such that $G-x$ is a tree. Furthemore, 
a graph $G$ is called an {\it $\ell$-apex tree} if there exists a vertex subset $A\subset V(G)$ of cardinality $\ell$ such that $G-A$ is a tree and there is no other subset of smaller cardinality with this property \cite{ap-xu,ap-kl}. 

In \cite{ap-kl} extremal values of (additively and multiplicatively) weighted Harary indices of apex and $\ell$-apex trees were studied. Extremal values of some other topological indices of $\ell$-apex trees were recently explored in \cite{hind21} and \cite{apex20}. In the later authors studied the generalized Wiener index
and derived the following result in which $K_{\ell}+T$
denotes the join of a complete graph $K_{\ell}$ and a tree $T$ on $n-\ell$ vertices.

\begin{theorem}
	\label{cor:min}
	Let $G$ be an $\ell$-apex tree on $n$ vertices, where $\ell\ge1$ and $n\ge \ell+2$,
	and let $\alpha\ne 0$. Then, the following two claims hold:
	\begin{itemize}
		\item
		If $\alpha>0$ then $W^{\alpha}(G)$ has the minimum value if and only if
		$G=K_{\ell}+T$, where $T$ is any tree on $n-\ell$ vertices;
		\item If $\alpha<0$ then $W^{\alpha}(G)$ has the maximum value if and only if
		$G=K_{\ell}+T$, where $T$ is any tree on $n-\ell$ vertices.
	\end{itemize}
	Moreover, in the extremal case
	$$W^{\alpha} (G)=(n^2-2n\ell-3n+\ell^2+3\ell+2)\,2^{\alpha-1}
		+\\(2n\ell+2n-\ell^2-3\ell-2)\,2^{-1}.$$
\end{theorem}

Observe that for $\alpha=1$ the invariant $W^{\alpha}$ is the Wiener index, and by Theorem~{\ref{cor:min}} the extremal value is
$$ W(G)=(2n^2-2n\ell-4n+\ell^2+3\ell+2)\,2^{-1}. $$

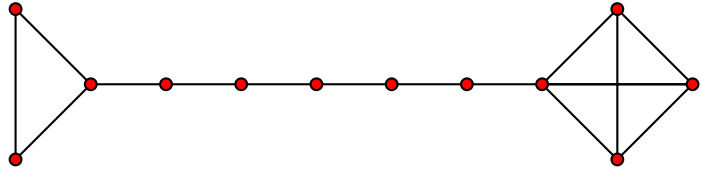
\begin{figure}[h]
\begin{center}
\begin{tikzpicture}[scale=1,style=thick]

		\node [My Style, name=a]   at (-1,1) {};		
		\node [My Style, name=b]   at (0,0) {};
		\node [My Style, name=c]   at (-1,-1) {};		
		\node [My Style, name=d]   at (1,0) {};
		\node [My Style, name=e]   at (2,0) {};
		\node [My Style, name=f]   at (3,0) {};
		\node [My Style, name=g]   at (4,0) {};
		\node [My Style, name=h]   at (5,0) {};
		\node [My Style, name=i]   at (6,0) {};
		\node [My Style, name=j]   at (7,1) {};
		\node [My Style, name=k]   at (8,0) {};
		\node [My Style, name=l]   at (7,-1) {};

		\draw[] (b)--(a)--(c)--(b)--(d)--(e)--(f)--(g)--(h)--(i)--(k)--(l);
		\draw[](k)--(j)--(l)--(i)--(k);
		\draw[](i)--(j);

\end{tikzpicture}
\end{center}
\caption{The graph $D_{5}(3,4)$.}

\label{dum}
\end{figure}

Recall that a dumbbell graph is a graph comprised of two disjoint cliques connected by a path. More precisely, a dumbbell graph $D_c(a,b)$ is a graph obtained from  a path $P_{c}=v_1v_2\cdots v_c$ and disjoint complete graphs $K_a$ and $K_b$ by connecting $v_1$ to a vertex of $K_a$ and connecting $v_c$ to a vertex of $K_b$, see Figure~\ref{dum} for $D_5(3,4)$. The order of so constructed graph is $a+b+c$.  Note that without loss of generality, we can always assume that $a,b\ne 2$.

\begin{theorem}
\label{cor:max}
Let $G$ be an apex tree on $n\ge 3$ vertices,
and let $\alpha\ne 0$.
\begin{itemize}
\item
If $\alpha>0$ then $W^{\alpha}(G)$ has the maximum value if and only if
$G=D_{n-4}(3,1)$;
\item
If $\alpha<0$ then $W^{\alpha}(G)$ has the minimum value if and only if
$G=D_{n-4}(3,1)$.
\end{itemize}
Moreover, in the extremal case
$$
W^{\alpha}(G)=1+\sum_{i=1}^{n-2}(n-i)i^{\alpha}.
$$
\end{theorem}


In~\cite{apex20} the following conjecture was proposed.

\begin{conjecture}
Let $G$ be an $\ell$-apex tree on $n$ vertices, where $\ell\ge3$ 
and $n\ge \ell+1$, such that $G$ has maximum Wiener index. Then $G$ is the balanced dumbbell graph, i.e. $G\cong D_c(a,b)$, 
where $a=\lceil \ell/2\rceil$, $b=\lfloor \ell/2\rfloor$, and $c=n-\ell$.
\end{conjecture}


%


\section{Wiener index of line graphs}

The line graph $L(G)$ of a graph $G$ is defined as a graph whose vertex set
coincides with the set of edges of $G$ and two vertices of $L(G)$ are adjacent
if and only if the corresponding edges are incident in $G$. Higher iterations
of the line graph are defined recursively.
\[
L^{k}(G)=\left\{
\begin{array}
[c]{ll}%
G & \text{for }k=0,\\
L(L^{k-1}(G)) & \text{for }k>0.
\end{array}
\right.
\]
Van Rooij and Wilf \cite{vanR} showed that for the sequence $$G,L(G),L(L(G)),L(L(L(G))),\ldots$$ only four options are possible. If $G$ is a cycle graph, then $L(G)$ and each subsequent graph in this sequence is isomorphic to $G$ itself. If $G$ is a claw $K_{1,3}$, then $L(G)=C_3$ and consequently the same holds for all subsequent graphs in the sequence. For a path we have $L(P_n)=P_{n-1}$, $L^2(P_n)=P_{n-2}$,
\ldots, $L^{n-1}(P_n)=P_1$ and $L^k(P_n)$ is an empty graph if $k\ge n$.
In all the remaining cases the order of the graphs in the sequence increases without bound.

\medskip

The following problem was proposed by Gutman~\cite{g-dlg-96}.
\begin{problem}
Find an $n$-vertex graph $G$ whose line graph $L(G)$ has maximum Wiener index.
\end{problem}

Supported by a result from \cite{dgms-ewig-09}, we pose the following conjecture (see also \cite{mathasp}).

\begin{conjecture}
Among all graphs $G$ on $n$ vertices, $W(L(G))$ attains maximum for some dumbbell graph on $n$ vertices.
\end{conjecture}

Similar conjecture was proposed for bipartite graphs \cite{mathasp}. Let us call a graph a {\em barbell} graph if it is comprised of two disjoint complete bipartite graphs connected by a path.

\begin{conjecture}
Let $n$ be large. Among all bipartite graphs $G$ on $n$ vertices, $W(L(G))$ attains maximum for some barbell graph on $n$ vertices.
\end{conjecture}

A related question we pose is the following.

\begin{problem} For given $n$ and $k$, find graphs $G$ on $n$ vertices with
the extremal value of $W(L^k(G))$.
\end{problem}

\medskip

Dobrynin and Mel'nikov~\cite{DM3} proposed to estimate the extremal values for the ratio $\frac{W(L^k(G))}{W(G)}$, for a graph $G$ on $n$ vertices and explicitly stated the case $k=1$ as a problem. The minimal value was given in \cite{miAMC}.

\begin{theorem}
Among all connected graphs on $n$ vertices, the fraction $\frac{W(L(G))}{W(G)}$ is
minimum for the star $S_n$, in which case $\frac{W(L(G))}{W(G)}=\frac{n-2}{2(n-1)}$.
\end{theorem}

The problem was recently solved also for the maximal value \cite{SSedlar21} .

\begin{theorem}
For a graph $G$ on $n$ vertices it holds that 
$\frac{W(L(G))}{W(G)}\leq\binom{n-1}{2}$
with equality if and only if $G=K_{n}.$
\end{theorem}

For $k>1$ the problem remains open.

\begin{problem}
\label{mel}
 Find  $n$-vertex graphs $G$ with extremal values of $\frac{W(L^k(G))}{W(G)}$ for $k\geq 2$.
\end{problem}

Note that the line graph of $K_n$ has the greatest number of vertices, and
restricting to bipartite graphs, the (almost) balanced complete bipartite graphs
have line graphs with most vertices, so $K_{\lfloor n/2\rfloor, \lceil n/2\rceil}$ could be the graph attaining maximal value in this class of graphs. 
It is expected that the minimum value should be attained by $P_n$, since this is the only graph whose line graph decreases in size, see a conjecture from \cite{mathasp}.

\begin{conjecture}
Let  $k\ge 2$ and let $n$ be large. Among all graphs $G$ on $n$ vertices, $\frac{W(L^k(G))}{W(G)}$ attains the maximum for $K_n$, and 
it attains the minimum for $P_n$.
\end{conjecture}

The above conjecture is supported by a result from \cite{HKS21}, where it was proved that among all trees on $n$ vertices the path $P_n$ has the smallest value of this ratio for $k \geq 3$, and it was conjectured that the same holds also in the case $k=2$. Another related problem is the following.

\begin{problem}
For various $\ell$ and $k$ find the extremal graphs for the ratio $\frac{W(L^k(G))}{W(L^l(G))}$.
\end{problem}

%
%
\section{Graphs with prescribed number of blocks }

A graph is {\em non-separable} if it is connected and has no
cut-vertices, i.e.~either it is $2$-connected or it is $K_2$.
A {\it block} of $G$ is a maximal non-separable subgraph of $G$.
As known, the $n$-path $P_n$, which has $n-1$ blocks, has the maximum Wiener index in the class of graphs on $n$ vertices, and among graphs on $n$ vertices that have just one block, the $n$-cycle has the largest Wiener index. The ordering of trees with respect to decreasing Wiener index is known up to the 17th maximum Wiener index \cite{1,9}, and the increasing ordering up to the 15th maximum Wiener index \cite{3}.

\begin{figure}[h]
\begin{center}
\begin{tikzpicture}[scale=0.45,style=thick]

		\node [My Style, name=a]   at (0,2) {};		
		\node [My Style, name=b]   at (1,4) {};
		\node [My Style, name=c]   at (3,5) {};		
		\node [My Style, name=d]   at (5,5) {};
		\node [My Style, name=e]   at (7,4) {};
		\node [My Style, name=f]   at (8,2) {};
		\node [My Style, name=g]   at (4,0) {};

		\draw[] (a)--(b)--(c);
		\draw[dotted] (c)--(d);
		\draw[] (d)--(e)--(f)--(g)--(a);
		\draw[] (f)--(a);

		\node [My Style, name=a1]   at (10,2) {};		
		\node [My Style, name=b1]   at (11,4) {};
		\node [My Style, name=c1]   at (13,5) {};		
		\node [My Style, name=d1]   at (15,5) {};
		\node [My Style, name=e1]   at (17,4) {};
		\node [My Style, name=f1]   at (18,2) {};
		\node [My Style, name=g1]   at (14,0) {};
		\node [My Style, name=h1]   at (14,2) {};	
	
		\draw[] (a1)--(b1)--(c1);
		\draw[dotted] (c1)--(d1);
		\draw[] (d1)--(e1)--(f1)--(g1)--(a1);
		\draw[] (f1)--(h1)--(a1);

		\node [My Style, name=a2]   at (20,2) {};		
		\node [My Style, name=b2]   at (21,4) {};
		\node [My Style, name=c2]   at (23,5) {};		
		\node [My Style, name=d2]   at (25,5) {};
		\node [My Style, name=e2]   at (27,4) {};
		\node [My Style, name=f2]   at (28,2) {};
		\node [My Style, name=g2]   at (24,0) {};
		\node [My Style, name=h2]   at (24,2) {};	
	
		\draw[] (a2)--(b2)--(c2);
		\draw[dotted] (c2)--(d2);
		\draw[] (d2)--(e2)--(f2)--(g2)--(a2);
		\draw[] (f2)--(h2)--(a2);
		\draw[] (h2)--(g2);
		
\end{tikzpicture}
\end{center}
\caption{Graphs $H_{n,1,2}$, $H_{n,2,2}$ and $H^+_{n,2,2}$.}

\label{hgraf}
\end{figure}
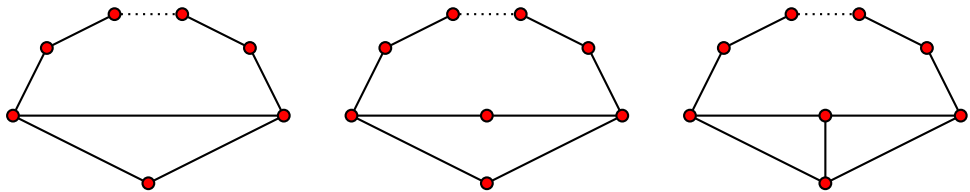

Bessy et al.~\cite{bess3} studied the ordering of $n$-vertex graphs with just one block (i.e.~$2$-vertex connected graphs)  with respect to decreasing Wiener index.
Let $1 \leq p \leq q \leq n - p -q + 1$ and $q > 1$. The notation $H_{n,p,q}$ stands for the graph on $n$ vertices comprised of
three internally disjoint paths with the same end-vertices, where the first path has length $p$, the second one has length $q$, and the last one has length $n-p-q+1$. Obviously $H_{n,1,2}$ is a graph
obtained from $C_n$ by introducing a new edge connecting two vertices at distance two on the cycle, and $H_{n,2,2}$ is a graph that is obtained from a $4$-cycle by connecting opposite vertices by a path of length $n-3$, see Figure~\ref{hgraf}.


In~\cite{bess3} it was shown that among graphs on $n$ vertices that have just one block, $H_{n,1,2}$ has the second largest Wiener index if $n \neq 6$. If $n \geq 11$, the third extremal graph is $H_{n,2,2}$. 
The authors also give conjectures on the graphs with 4th and 5th greatest Wiener index in the class of $2$-connected graphs.
Let $H^+_{n,2,2}$ be the graph obtained from $H_{n,2,2}$ by inserting an edge between two vertices that are at distance $1$ from the
vertices of degree $3$, see the third graph in Figure~\ref{hgraf}. Then $H^+_{n,2,2}$ has Wiener index exactly $1$ less than $H_{n,2,2}$, so it is the fourth $2$-connected graph by decreasing Wiener index for $n = 9$ and $n \geq 11$, but it may not be  unique. However, the following can be true.

\begin{conjecture} For $n$ large
enough, $H^+_{n,2,2}$ is the graph with the 4th largest Wiener index among blocks on $n$ vertices.
\end{conjecture}

\begin{conjecture}
For $n$ large
enough, $H_{n,1,3}$ is the graph with the 5th largest Wiener index among blocks on $n$ vertices.
\end{conjecture}

Bessy et al.~\cite{bess1} studied a general problem of finding the maximum possible value of Wiener index among graphs on $n$ vertices with fixed number of blocks. They showed that among all graphs on $n$ vertices which have $p\ge 2$ blocks, the maximum Wiener index is attained by a graph comprised of two cycles
joined by a path, where one or both cycles can be replaced 
by a single edge. To be more specific, we need the following notation.

If $G$ is a connected graph and $v$ is a cut-vertex that partitions
$G$ into subgraphs $G_1$ and $G_2$, i.e., $G=G_1\cup G_2$ and
$G_1\cap G_2 = \{v\}$, then we write $G=G_1\circ_v G_2$. For simplicity reasons, by $C_2$ we mean the complete graph $K_2$.

\begin{theorem}
\label{thm:main}
Let $n$ and $p$ be numbers such that $n>p>1$.
Among all graphs on $n$ vertices with $p$ blocks, the maximum
Wiener index is attained by the graph
$C_a\circ_u P_{p-1}\circ_v C_b$
for some integers $a \ge 2$ and $b \ge 2$, where $a+b=n-p+3$,
and $u$ and $v$ are distinct endvertices of $P_{p-1}$.
\end{theorem}

Note that $C_a$ or $C_b$ can also be edges, and then we
obtain $C_{n-p+1}\circ_u P_{p}$, which is a graph composed of one cycle
with an attached path, or $P_n$ if both $C_a$ and $C_b$ are
edges.

In \cite{bess2} the authors provide further details  by determining the sizes of
$a$ and $b$ in the extremal graphs for each $n$ and $p$. Roughly speaking, if $n$ is bigger than $5p-7$, then the extremal graphs is obtained for  $a=2$, i.e.~the graph is a path glued to a cycle. For values $n=5p-8$ and $5p-7$, there is more than one extremal graph. And when $n<5p-8$, the extremal graphs is again unique with $a$ and $b$ being equal or almost equal depending on the congruence of $n-p$ modulo $4$. 

\bigskip

\bigskip\noindent\textbf{Acknowledgments.}~~The first author
acknowledges
partial support by Slovak research grants VEGA 1/0567/22, VEGA
1/0206/20,
APVV--19--0308, APVV--17--0428. All authors acknowledge partial support
of the Slovenian research agency
ARRS program\ P1-0383 and ARRS project J1-1692.


\begin{thebibliography}{99}

\bibitem{akh}
M.~Akhmejanova, K.~Olmezov, A.~Volostnov, I.~Vorobyev, K.~Vorob'ev, Y.~Yarovikov,
Wiener index and graphs, almost half of whose vertices satisfy \v{S}olt\'{e}s property, 
{\it Appl. Math Comput.\/} {\bf 325} (2023) 37--42.



\bibitem{hind21}
A.~Ali, W.~Iqbal, Z.~Raza, E.~E.~Ali, J.~Liu, F.~Ahmad, Q.~A.~Chaudhry, 
Some Vertex/Edge-Degree-Based Topological Indices of $r$-Apex Trees, 
Journal of Mathematics, vol. 2021, Article ID 4349074, 8 pages, 2021. 
https://doi.org/10.1155/2021/4349074

\bibitem{win-dimension} 
Y.~Alizadeh, V.~Andova, S.~Klav\v{z}ar, R. \v{S}krekovski,
Wiener Dimension: Fundamental Properties and (5,0)-Nanotubical Fullerenes,
MATCH Commun. Math. Comput. Chem 72 (2014) 279--294.

\bibitem{AloDankel}
A.~Alochukwu, P.~Dankelmann,
Wiener index in graphs with given minimum degree and maximum degree,
Discrete Mathematics \& Theoretical Computer Science, vol. 23 no. 1 (2021).






\bibitem{basic}
N.~Ba\v{s}i\'{c}, M.~Knor, R.~\v{S}krekovski,
On regular graphs with Šoltés vertices,
manuscript.

\bibitem{bess1}
S.~Bessy, F.~Dross, K.~Hri\v{n}\'{a}kov\'{a}, M.~Knor, R.~\v{S}krekovski, 
The structure of graphs with given number of blocks and the maximum Wiener index, 
{\it J. of Combin. Optim.\/} {\bf 39} (2020) 170--184.


\bibitem{bess2}
S.~Bessy, F.~Dross, K.~Hri\v{n}\'{a}kov\'{a}, M.~Knor, R.~\v{S}krekovski, 
Maximal Wiener index for graphs with prescribed number of blocks, 
{\it Appl. Math Comput.\/} {\bf 380} (2020) 125274.


\bibitem{bess3}
S.~Bessy, F.~Dross, M.~Knor, R.~\v{S}krekovski, 
Graphs with the second and third maximum Wiener index over the $2$-vertex connected graphs,
{\it Discrete Appl. Math.\/} {\bf 284} (2020) 195--200.

\bibitem{Bok1}
J.~Bok, N.~Jedličkov\'{a}, J.~Maxov\'{a},
On relaxed \v{S}olt\'{e}s's problem,
{\it Acta Math. Univ. Comenianae\/} Vol. LXXXVIII, 3 (2019) 475--480.

\bibitem{Bok2}
J.~Bok, N.~Jedličkov\'{a}, J.~Maxov\'{a},
A relaxed version of \v{S}olt\'{e}s's problem and cactus graphs,
Bull. Malays. Math. Sci. Soc. 44 (2021) 3733--3745. 


\bibitem{bonc1}
D. Bonchev,
{\em On the complexity of directed biological networks},
SAR QSAR Environ Res. \textbf{14} (2003), 199--214.

\bibitem{bonc2}
D. Bonchev,
{\em Complexity of Protein-Protein Interaction Networks, Complexes and Pathways},
in Handbook of Proteomics Methods, M. Conn, ed. Humana, New York, (2003), 451--462.

\bibitem{bozo}
V.~Bo\v{z}ovi\'{c}, \v{Z}.~Kovijani\'{c} Vukićevi\'{c}, G.~Popivoda,
R.~Pan, X.~Zhang,
Extreme Wiener indices of trees with given number of vertices
of maximum degree,
\emph{Discrete Appl. Math.} \textbf{304} (2021) 23--31.


\bibitem{chen1}
Y.~Chen, X.~Lin, X.~Zhang,
The extremal average distance of cubic graphs.
{\em Journal of Graph Theory}, manuscript.


\bibitem{dgms-ewig-09}
P.~Dankelmann, I.~Gutman, S.~Mukwembi, H.C.~Swart,
The edge-Wiener index of a graph,
\emph{Discrete Math.} \textbf{309} (2009), 3452--3457.

\bibitem{das17}
K.~C.~Das, M.~J.~Nadjafi-Arani,
On maximum Wiener index of trees and graphs with given radius,
J Comb Optim 34 (2017)  574--587.

\bibitem{CambieDiam}
S.~Cambie,
An asymptotic resolution of a problem of Plesn\'{i}k,
J. Comb. Theory. Ser. B, 145 (2020) 341--358.

\bibitem{Cambie21rad}
S.~Cambie,
Extremal total distance of graphs of given radius I,
J. Graph Theory. 97 (2021) 104--122.

\bibitem{CambHas}
S.~Cambie, J.~Haslegrave, 
On the relationship between variable Wiener index and variable
Szeged index,
\emph{Appl. Math. Comput.} 431 (2022) 127320.

\bibitem{Chen}
Y.~Chen, B.~Wu, X.~An,
Wiener Index of Graphs with Radius Two,
\emph{ISRN Combinatorics}, Article ID 906756 (2013), 5 pages.

\bibitem{dankelmann22}
P.~Dankelmann,
On the Wiener Index of Orientations of Graphs,
manuscript, arXiv:2209.08946 [math.CO].

\bibitem{DeLa}
E.~DeLaVi\~{n}a, B.~Waller,
Spanning trees with many leaves and average distance, 
\emph{Electronic J. Combin.} \textbf{15} (2014), R33 p.16.

\bibitem{1}
H.~K.~Deng, 
The trees on $n\geq 9$ vertices with the first to seventeenth greatest Wiener indices are chemical graphs, 
MATCH Commun. Math. Comput. Chem. 57 (2007) 393--402.

\bibitem{surv1}
A.~A.~Dobrynin, R.~Entringer, I.~Gutman,
Wiener index of trees: Theory and applications,
\emph{Acta Appl. Math.} \textbf{66}(3) (2001), 211--249.

\bibitem {DobGut95}
A.~Dobrynin, I.~Gutman,
Solving a problem connected with distances in graphs,
\emph{Graph Theory Notes N. Y.}, \textbf{28} (1995), 21--23.

\bibitem{petra}
A.~A.~Dobrynin, I.~Gutman, S.~Klav\v zar, P. \v Zigert, 
Wiener index of hexagonal systems,
\emph{Acta Appl. Math.} \textbf{72} (2002), 247--294.


\bibitem{DM3}
A.~A.~Dobrynin, L.~S.~Mel'nikov,
Wiener index of line graphs, in I. Gutman, B. Furtula (Eds.)
{\em Distance in Molecular Graphs – Theory}, Univ. Kragujevac, Kragujevac (2012), 85--121.

\bibitem{DongZhou} 
H.~Dong, B.~Zhou,
Maximum Wiener index of unicyclic graphs with fixed maximum degree,
{\em Ars Combin.\/} {\bf  103}  (2012) 407--416.

\bibitem{3} 
H.~Dong, X.~Guo, 
Ordering trees by their Wiener indices, 
MATCH Commun. Math. Comput. Chem. 56 (2006) 527--540.

 
\bibitem{n1}
R.~C.~Entringer, D.~E.~Jackson,  D.~A.~Snyder, Distance in graphs, 
{\it Czechoslovak Math. J.\/} {\bf 26} (1976) 283--296.

\bibitem{EJ}
G.~Exoo, R.~Jajcay,
Dynamic cage survey,
{\it Electronic J. Combin. Dynamic Survey DS16\/} (2013).

\bibitem{counterex}
Y.~Fang, Y.~Gao,
Counterexamples to the conjecture on orientations of graphs
with minimum Wiener index,
\emph{Discrete Appl. Math.} \textbf{232} (2017) 213--220.

\bibitem{fish} 
M.~Fischermann, A.~Hoffmann, D.~Rautenbach, L.~Sz\'{e}kely, L.~Volkmann, 
Wiener index versus maximum degree in trees, 
\emph{Discrete Appl. Math.} \textbf{122} (2002) 127--137.

\bibitem{n2}
I.~Gutman, A property of the Wiener number and its modifications, {\it Indian J. Chem.\/} {\bf 36A} (1997) 128--132.

\bibitem{Gut94}
I. Gutman,
A formula for the Wiener number of trees and its extension to graphs containing cycles,
\emph{Graph Theory Notes N. Y.} \textbf{27} (1994), 9--15.

\bibitem{g-dlg-96}
I.~Gutman,
Distance of line graphs,
\emph{Graph Theory Notes N. Y.} \textbf{31} (1996), 49--52.



\bibitem{n3}
I.~Gutman,  W.~Linert, I.~Lukovits,  A.~A.~Dobrynin, Trees with extremal hyper-Wiener index: Mathematical basis and chemical applications, {\it J. Chem. Inf. Comput. Sci.\/} {\bf 37} (1997) 349--354.

\bibitem{GutVukZer04}
I. Gutman, D. Vuki\v cevi\'c, J. \v Zerovnik,
A class of modified Wiener indices,
\emph{Croat. Chem. Acta} \textbf{77} (2004), 103--109.

\bibitem{henol}
M.~A.~Henning, O.~R.~Oellermann,
The average connectivity of a digraph. 
\emph{Discrete Appl. Math.} \textbf{140} (2004), 143--153.

\bibitem{Hosoya}
H.~Hosoya,
Topological index. A newly proposed quantity characterizing the topological
nature of structural isomers of saturated hydrocarbons,
Bull.~Chem.~Soc.~Jpn. 44(9) (1971) 2332--2339. 

\bibitem{HKS19}
K.~Hri\v{n}\'{a}kov\'{a}, M.~Knor, R.~\v{S}krekovski,
An inequality between variable Wiener index and variable Szeged index,
{\it Appl. Math. Comput.\/} {\bf 362} (2019) 124557.

\bibitem{HKS21}
K.~Hri\v{n}\'{a}kov\'{a}, M.~Knor, R.~\v{S}krekovski,
On a conjecture about the ratio of Wiener index in iterated line graphs,
{\it The Art of Discrete and Applied Mathematics\/} {\bf 1} (2018),
https://doi.org/10.26493/2590-9770.1257.dda

\bibitem{Hu21}
Y.~Hu, Z.~Zhu, P.~Wu, Z.~Shao, A.~Fahad,
On investigations of graphs preserving the Wiener index upon vertex removal,
AIMS Mathematics, \textbf{6(12)} (2021) 12976–12985.


\bibitem{jel}
F.~Jelen, 
Superdominance order and distance of trees, Doctoral thesis, RWTH
Aachen, Germany, 2002.

\bibitem{tri} 
F.~Jelen, E. Trisch, 
Superdominance order and distance of trees with bounded maximum degree,
\emph{Discrete Appl. Math.} \textbf{125} (2003), 225--233.

\bibitem{KhaDesKal95}
P.~V.~Khadikar, N.~V.~Deshpande, P.~P.~Kale, A.~Dobrynin, I.~Gutman, G.~D\"om\"ot\"or,
The Szeged index and an analogy with the Wiener index,
\emph{J. Chem. Inf. Comput. Sci.} \textbf{35 (3)} (1995), 547--550.

\bibitem{KlaRajGut96}
S.~Klav\v zar, A.~Rajapakse, and I.~Gutman,
The Szeged and the Wiener index of graphs,
\emph{Appl. Math. Lett.} \textbf{9 (5)} (1996), 45--49.

\bibitem{apex20}
M.~Knor, M.~Imran, M.~K.~Jamil, R~\v{S}krekovski,
Remarks on distance based topological indices for $\ell$-apex trees,
{\it Symmetry\/} {\bf 12(5):802} (2020). 

\bibitem{Majst-1}
M.~Knor, S.~Majstorovi\'{c}, R~\v{S}krekovski,
Graphs whose Wiener index does not change when a speciffic vertex is removed,
{\it Discrete Appl. Math.\/} {\bf 238} (2018) 126--132.


\bibitem{Majst-2}
M.~Knor, S.~Majstorovi\'{c}, R.~\v{S}krekovski,
Graphs preserving Wiener index upon vertex removal,
{\it Appl. Math. Comput.\/} {\bf 338} (2018) 25--32.

\bibitem{Majst-3}
M.~Knor, S.~Majstorovi\'{c}, R.~\v{S}krekovski,
Some results on wiener index of a graph: an overview,
2nd Croatian Combinatorial Days, page 49.

\bibitem{KS_chapter}
M.~Knor, R.~\v Skrekovski, Wiener index of line graphs,  in  M. Dehmer and F. Emmert-Streib (Eds.),
{\em Quantitative Graph Theory: Mathematical Foundations and Applications)}, CRC Press (2014), 279--301.

\bibitem{KSgrid} 
M.~Knor, R.~\v Skrekovski, On maximum Wiener index of directed grids,
https://doi.org/10.48550/arXiv.2201.11958.

\bibitem{mathasp} 
M.~Knor, R.~\v Skrekovski, A.~Tepeh,
Mathematical aspects of Wiener index, 
{\it Ars Math. Contemp.\/} {\bf 11} (2016) 327--352.

\bibitem{KST}
M.~Knor, R.~{\v S}krekovski, A.~Tepeh,
Orientations of graphs with maximum Wiener index, 
Discrete Appl. Math. 211 (2016) 121--129. 

\bibitem{KST2}
M.~Knor, R.~{\v S}krekovski, A.~Tepeh, 
Some remarks on the Wiener index of digraphs, 
Appl. Math. Comput. 273 (2016) 631--636.

\bibitem{KST3}
M.~Knor, R.~{\v S}krekovski, A.~Tepeh, 
Digraphs with large maximum Wiener index,
Appl. Math. Comput. 284 (2016), 260--267.


\bibitem{KST-dir-survey} 
M.~Knor, R.~{\v S}krekovski, A.~Tepeh,
Wiener index of digraphs. In: Gutman, I.,  Furtula, B.,   Das,  K. C., Milovanovic, E.,  Milovanovic, I. (eds.)  Bounds in Chemical Graph Theory -- Advances, pp. 141--153. Univ. Kragujevac, Kragujevac (2017).

\bibitem{miAMC}
M.~Knor, R.~{\v S}krekovski, A.~Tepeh, 
An inequality between the edge-Wiener index and the Wiener index of a graph,
\emph{Appl. Math. Comput.} \textbf{269} (2015), 714--721.


\bibitem{mi-chem}
M.~Knor, R.~\v Skrekovski, A.~Tepeh,
Chemical graphs with the minimum value of Wiener index, 
{\it MATCH Commun. Math. Comput. Chem.\/} {\bf 81} (2019) 119--132.

\bibitem{ladder21}
T.~Kraner \v{S}umenjak, S.~\v{S}pacapan,  D. \v{S}tesl,
A proof of a conjecture on maximum Wiener index of oriented ladder  graphs,
J. Appl. Math. Comput. (2021). https://doi.org/10.1007/s12190-021-01498-w.

\bibitem{ZanaBul} 
\v{Z}.~Kovijani\'{c} Vuki\'{c}evi\'{c}, L.~Bulatovi\'{c},
On the variable Wiener - Szeged inequality,
Discrete Appl. Math. 307 (2022) 15--18. 


\bibitem{LiWu}
Z.~Li, B.~Wu,
Orientations of graphs with maximum Wiener index,
manuscript.


\bibitem{Lin14}
H.~Lin,
A note on the maximal Wiener index of trees with given number of vertices of
maximum degree,
{\it MATCH Commun. Math. Comput. Chem.\/} {\bf 72} (2014) 783--790.

\bibitem{liupan}
H.~Liu, X.~Pan,
On the Wiener index of trees with fixed diameter,
{\it MATCH Commun. Math. Comput. Chem.\/} {\bf 60} (2008) 85--94.

\bibitem{9}
M.~Liu, B.~Liu, Q.~Li, 
Erratum to 'The trees on $n \geq 9$ vertices with the first to seventeenth greatest Wiener indices are chemical graphs', 
MATCH Commun. Math. Comput. Chem. 64 (2010) 743--756.


\bibitem{dW}
E.~Loz, H.~P{\'e}res-Ros{\'e}s, G.~Pineda-Villavicencio, The degree diameter problem for general graphs,
{\it Combinatorics Wiki\/}, 24 May 2017, 15:30 UTC,
$\langle$http://combinatoricswiki.org/wiki/The\_{}Degree\_{}Diameter\_{}Problem\_{}for\_{}General\_{}Graphs$\rangle$
[accessed 3 February 2018]


\bibitem{MS}
M.~Miller, J.~{\v S}ir{\'a}{\v n},
Moore graphs and beyond: A survey of the degree-diameter problem,
{\it Electronic J. Combin. Dynamic Survey D14\/} (2005).

\bibitem{muk}
S.~Mukwembi, T.~Vert{\'i}k, 
Wiener index of trees of given order and diameter at most $6$,
\emph{Bull. Aust. Math. Soc.} \textbf{89} (2014), 379--396.


\bibitem{moon}
J. W. Moon,  On the total distance between nodes in
tournaments, Discrete Math. \textbf{151} (1996), 169--174.



\bibitem{P}
J.~Plesn{\'\i}k,
On the sum of all distances in a graph or digraph,
\emph{J. Graph Theory} \textbf{8} (1984), 1--21.


\bibitem{SSedlar21}
J.~Sedlar, R.~\v Skrekovski,
A note on the maximum value of $W(L(G)) / W(G)$, manuscript.

\bibitem{sun19}
Q.~Sun, B.~Ikica, R.~\v{S}krekovski, V.~Vuka\v{s}inovi\'{c},
Graphs with a given diameter that maximize the Wiener index,
{\it Appl. Math. Comput.\/} {\bf 356} (2019) 438--448.

\bibitem{soltes}
L'.~\v{S}olt\'{e}s,
Transmission in graphs: A bound and vertex removing,
{\it Math. Slovaca \/} {\bf 41} (1991) 11--16.

\bibitem{spiro21}
S.~Spiro,
The Wiener Index of signed graphs,
{\it Appl. Math. comput.\/} {\bf 416} (2022) 126755.

\bibitem{stev}
D.~Stevanovi\'{c}, 
Maximizing Wiener index of graphs with fixed maximum degree,
\emph{MATCH Commun. Math. Comput. Chem.} \textbf{60} (2008), 71--83.

\bibitem{stevan20}
D.~Stevanovi\'{c}, N.~Milosavljevi\'{c}, D.~Vuki\v{c}evi\'{c},
A few examples and counterexamples in spectral graph theory,
Discussiones Mathematicae Graph Theory 40 (2020) 637--662.

\bibitem{S}
Stojmenovi\'c, I. 
Honeycomb networks: Topological properties
and communication Algorithms,
{\em IEEE Trans. Parallel Distrib. Syst.} {\bf 1997} {\em 8} 1036--1042.

\bibitem{vanR}
A.~C.~M. van Rooij, H.~S.~Wilf, 
The interchange graph of a finite graph,
{\em Acta Mathematica Academiae, Scientiarum Hungaricae} {\bf 16} (1965) 263--269.

\bibitem{n7}
H.~Wang, The extremal values of the Wiener index of a tree with given degree sequence,
{\it Discrete Appl. Math.\/} {\bf 156} (2009) 2647--2654.


\bibitem{wang2}
S.~Wang, X.~Guo,
Trees with extremal Wiener indices,
\emph{MATCH Commun. Math. Comput. Chem.} \textbf{60} (2008) 609--622.

\bibitem{Wiener}
H.~Wiener,
Structural determination of paraffin boiling points,
\emph{J. Amer. Chem. Soc.} \textbf{69} (1947) 17--20.

\bibitem{You}
Z.~You, B.~Liu, 
Note on the minimal Wiener index of connected graphs with $n$ vertices and radius $r$,
\emph{MATCH Commun. Math. Comput. Chem.} \textbf{66} (2011), 343--344.

\bibitem{Gut-sur}
K.~Xu, M.~Liu, K.C.~Das, I.~Gutman, B.~Furtula, 
A survey on graphs extremal with respect to distance-based topological indices
{\em MATCH Commun. Math. Comput. Chem.} \textbf{71} (2014), 461--508.

\bibitem{ap-xu}
K.~Xu, Z.~Zheng, K.~Ch.~Das, 
Extremal t-apex trees with respect to matching energy,
{\it Complexity\/} {\bf 21} (2015) 238--247.


\bibitem{ap-kl}
K.~Xu, J.~Wang, K.~Ch.~Das, S.~Klav\v{z}ar,
Weighted Harary indices of apex trees and $k$-apex trees,
{\it Discrete Appl. Math.\/} {\bf 189} (2015) 30--40.

\bibitem{n8} 
X.~D.~Zhang, Q.~Y.~Xing, L.~Q.~Xu and R.~Y.~Pang, The Wiener index of trees with given
degree sequences, {\it MATCH Commun. Math. Comput. Chem.\/} {\bf 60} (2008) 623--644.

\bibitem{zhang2} 
X.~D.~Zhang, Y.~Liu, M.~X.~Han,
Maximum Wiener index of trees with given degree sequence, 
\emph{MATCH Commun. Math. Comput. Chem.} \textbf{64} (2010) 661--682.


\end{thebibliography}
\end{document}